\documentclass{article}
\usepackage{amsfonts,amsmath,amssymb}
\usepackage{graphicx,epsfig}
\usepackage{cite,subfigure}
\usepackage{extarrows}
\usepackage{longtable,geometry}
\begin{document}


\newtheorem{theorem}{Theorem}[section]
\newtheorem{proposition}{Proposition}[section]
\newtheorem{lemma}{Lemma}[section]
\newtheorem{corollary}{Corollary}[section]
\newtheorem{remark}{Remark}[section]
\newtheorem{proof}{Proof:}

\renewcommand{\thesection}{\arabic{section}}
\renewcommand{\theequation}{\thesection.\arabic{equation}}
\renewcommand{\thetheorem}{\thesection.\arabic{theorem}}
\numberwithin{equation}{section}
\numberwithin{theorem}{section}
\numberwithin{proposition}{section}
\numberwithin{lemma}{section}
\numberwithin{remark}{section}
\setcounter{secnumdepth}{5}


\newcommand{\cl}{\centerline}
\newcommand{\sms}{\smallskip}
\newcommand{\ms}{\medskip}
\newcommand{\bs}{\bigskip}
\newcommand{\noi}{\noindent}
\newcommand{\itl}[1]{\textit{#1}}
\newcommand{\blf}[1]{\textbf{#1}}
\newcommand{\dsty}{\displaystyle}
\newcommand{\txty}{\textstyle}
\newcommand{\ssty}{\scriptstyle}
\newcommand{\tty}{\texttt}


\newcommand\Par{\mathhexbox278\,}


\newcommand{\al}{\alpha}
\newcommand{\Al}{\Alpha}
\newcommand{\be}{\beta}
\newcommand{\Be}{\Beta}
\newcommand{\Gm}{\Gamma}
\newcommand{\gm}{\gamma}
\newcommand{\dl}{\delta}
\newcommand{\Dl}{\Delta}
\newcommand{\lm}{\lambda}
\newcommand{\Lm}{\Lambda}
\newcommand{\kp}{\kappa}
\newcommand{\varep}{\varepsilon}
\newcommand{\vp}{\varphi}
\newcommand{\sig}{\sigma}
\newcommand{\Sig}{\Sigma}
\newcommand{\om}{\omega}
\newcommand{\Om}{\Omega}
\newcommand{\uom}{\mbox{\boldmath$\omega$}}
\newcommand{\btau}{\mbox{\boldmath$\tau$}}
\newcommand{\bnu}{\mbox{\boldmath$\nu$}}
\newcommand{\up}{\upsilon}
\newcommand{\z}{\zeta}


\newcommand{\df}[1]{\buildrel\mbox{\small def}\over{#1}}
\newcommand{\op}[1]{\buildrel\mbox{\tiny o}\over{#1}}
\newcommand{\db}{\prime\prime}
\newcommand{\bsl}{\backslash}
\newcommand{\lb}{\lbrack\!\lbrack}
\newcommand{\rb}{\rbrack\!\rbrack}
\newcommand\la{\langle}
\newcommand\ra{\rangle}
\newcommand{\ev}{\equiv}
\newcommand{\nev}{\not\equiv}
\newcommand{\nn}{\mathbb{N}}
\newcommand{\qq}{\mathbb{Q}}
\newcommand{\zz}{\mathbb{Z}}
\newcommand{\rr}{\mathbb{R}}
\newcommand{\rn}{\rr^N}
\newcommand{\cc}{\mathbb{C}}
\newcommand{\id}{\mathbb{I}}
\newcommand{\bo}{\mathbb{O}}

\newcommand{\amsb}[1]{\mathbb{#1}}
\newcommand{\mcl}[1]{\mathcal{#1}}
\newcommand{\bl}[1]{\mathbf{#1}}
\newcommand{\ov}[1]{\overline{#1}}
\newcommand{\wt}[1]{\widetilde{#1}}
\newcommand{\wh}[1]{\widehat{#1}}

\newcommand{\lra}{\longrightarrow}
\newcommand{\LLR}{\Longleftrightarrow}
\newcommand{\LRA}{\Longrightarrow}
\newcommand{\LLA}{\Longleftarrow}


\newcommand{\bbox}{\vrule height.6em width.6em 
depth0em} 
\newcommand{\os}{\vbox{\hrule \hbox{\vrule 
height.6em depth0pt 
\hskip.6em \vrule height.6em depth0em}
\hrule}} 


\newcommand{\dvg}{\operatorname{div}}
\newcommand{\curl}{\operatorname{curl}}
\newcommand{\supp}{\operatorname{supp}}
\newcommand{\essup}{\operatornamewithlimits{ess\,sup}}
\newcommand{\essinf}{\operatornamewithlimits{ess\,inf}}
\newcommand{\essosc}{\operatornamewithlimits{ess\,osc}}
\newcommand{\osc}{\operatornamewithlimits{osc}}
\newcommand{\sign}{\operatorname{sign}}
\newcommand{\loc}{\operatorname{loc}}
\newcommand{\diam}{\operatorname{diam}}
\newcommand{\dist}{\operatorname{dist}}
\newcommand{\card}{\operatorname{card}}
\newcommand{\meas}{\operatorname{meas}}
\newcommand{\spn}{\operatorname{span}}
\newcommand{\dtm}{\operatorname{det}}
%


\newcommand{\overlim}{\mathop{\overline{\lim}}\limits}
\newcommand{\underlim}{\mathop{\underline{\lim}}\limits}
\newcommand{\ttop}[2]{\genfrac{}{}{0pt}{}{#1}{#2}}
\newcommand{\bcu}{\mathop{\txty{\bigcup}}\limits}
\newcommand{\bca}{\mathop{\txty{\bigcap}}\limits}
\newcommand{\bsu}{\mathop{\txty{\sum}}\limits}
\newcommand{\pro}{\mathop{\txty{\prod}}\limits}


\newcommand{\pl}{\partial}
\newcommand{\ptt}{\frac{\pl}{\pl t}}
\newcommand{\ppx}{\frac\pl{\pl x}}
\newcommand{\dds}{\frac d{ds}}
\newcommand{\ddt}{\frac d{dt}}


\newcommand{\intl}{\int\limits}
\newcommand{\iintl}{\iint\limits}


\def\Xint#1{\mathchoice
    {\XXint\displaystyle\textstyle{#1}}%
    {\XXint\textstyle\scriptstyle{#1}}%
    {\XXint\scriptstyle\scriptscriptstyle{#1}}%
    {\XXint\scriptscriptstyle\scriptscriptstyle{#1}}%
    \!\int}
\def\XXint#1#2#3{\setbox0=\hbox{$#1{#2#3}{\int}$}
    \vcenter{\hbox{$#2#3$}}\kern-0.5\wd0}
\def\bint{\Xint-}
\def\dashint{\Xint{\raise4pt\hbox to7pt{\hrulefill}}}


\newcommand{\ovl}[3]{\int_{#1}^{#2}\kern-#3pt\raise4pt\hbox to7pt{\hrulefill}\ }

\newcommand{\ovll}[3]{\intl_{#1}^{#2}\kern-#3pt\raise4pt\hbox to7pt{\hrulefill}\ }

\newcommand{\tvl}[2]{\iint_{#1}\kern-#2pt\raise4pt\hbox to7pt{\hrulefill}\ }



\newcommand{\omt}{\Om_T}
\newcommand{\plo}{\partial\Omega}
\newcommand{\ovo}{\bar{\Om} }

%
\newcommand{\ci}[1]{C^\infty\!\left({#1}\right)}
\newcommand{\cio}[1]{C_o^\infty\!\left({#1}\right)}
\newcommand{\lloc}[1]{L_{\loc}\!\left({#1}\right)}
\newcommand{\xy}{|x-y|}


\newcommand{\intom}{\intl_{\Om}}
\newcommand{\intbo}{\intl_{\plo}}
\newcommand{\inom}{\int_{\Om}}
\newcommand{\inbo}{\int_{\plo}}
\newcommand{\intrn}{\intl_{\rn}}


\newcommand{\bye}{\end{document}}



%
%
%

%
%
%
%
%

%
%

\input harnack.mac
\begin{center}  {\huge\textbf{The Wiener Test for the Removability of the Logarithmic Singularity 
for the Elliptic PDEs with Measurable
Coefficients and Its Consequences }}
\par \medskip\bigskip\end{center}
\begin{center} {\Large\textsc{Ugur G. Abdulla}}
\par \medskip\bigskip\end{center}
\begin{center} {\large\noindent \textsc{Department of Mathematics, Florida Institute of Technology, Melbourne, Florida 32901}}
\par \medskip\bigskip\end{center}
{\bf Abstract.} This paper introduces the notion of $log$-regularity (or $log$-irregularity)
of the boundary point $\zeta$ (possibly $\zeta=\infty$) of the arbitrary open subset $\Omega$ of the Greenian deleted neigborhood of $\zeta$ in $\rr^2$ concerning second order uniformly elliptic equations with bounded and measurable coefficients, according as whether the  $log$-harmonic measure of $\zeta$ is null (or positive). A necessary and sufficient condition for the removability of the logarithmic singularity, that is to say for the existence of a unique solution
to the Dirichlet problem in $\Omega$ in a class $O(\log |\cdot - \zeta|)$ is established
in terms of the Wiener test for the $log$-regularity of $\zeta$.    
From a topological point of view, the Wiener test at $\zeta$ presents the minimal thinness criteria
of sets near $\zeta$ in minimal fine topology. Precisely, the open set $\Omega$ is a deleted neigborhood
of $\zeta$ in minimal fine topology if and only if $\zeta$ is $log$-irregular. From the probabilistic point of view, the Wiener test presents
asymptotic law for the $log$-Brownian motion near $\zeta$ conditioned on the logarithmic kernel with pole at $\zeta$.



\newpage
\section{Description of Main Results}\label{description of results,historical remarks}
\subsection{Introduction and Motivation}\label{E:1:1}
{\large
This paper introduces the notion of $log$-regularity of the boundary point 
$\zeta$ (possibly $\zeta=\infty$) and establishes a necessary and sufficient 
condition for the removability of the logarithmic singularity, or equivalently
for the unique solvability of the Dirichlet problem (DP)
in an arbitrary open subset $\Omega$ of the Greenian deleted neigborhood of $\zeta$
in $\rr^{2}$ for the uniformly elliptic 
equations in divergence form
\begin{equation}\label{Eq:W:1:1}
{\cal A}u=-(a_{ij}(x)u_{x_i})_{x_j}=0
\end{equation}
when the coefficients and boundary values 
are only supposed to be bounded and measurable. The particular case of the equation 
\eqref{Eq:W:1:1} is the Laplace equation, and the result of this paper 
is new for the classical case of harmonic functions.

In a recent paper \cite{Abdulla4} (see also \cite{Abdulla1})
 a notion of regularity (or irregularity)
of the point at infinity ($\infty$) is introduced for the unbounded open set $\Omega\subset\rr^{N}, N\ge 3$ concerning second order uniformly elliptic equations with bounded and measurable coefficients, according as whether the  ${\cal A}$- harmonic measure of $\infty$ is zero (or positive). A necessary and sufficient condition for the existence of a unique bounded solution
to the Dirichlet problem in an arbitrary open set of $\rr^{N}, N\ge 3$ is established
in terms of the Wiener test for the regularity of $\infty$. It coincides with the Wiener 
test for the regularity of $\infty$ in the case of the Laplace equation.  
This approach is not directly applicable to the case $N=2$. The Dirichlet problem 
in any Greenian open subset of $\rr^2$ always has a unique bounded solution independetly of the geometry 
of the boundary at $\infty$. The goal of this paper is to shed light upon the 
two-dimensional case and to demonstrate that the geometry of the boundary 
at $\infty$ has a similar crucial effect on the uniqueness of the solution to the Dirichlet problem
in a class of unbounded solutions with possibly a logarithmic growth rate at $\infty$.

Let
\[ D_1=\{x\in \rr^2: |x|>1 \}\]
If ${\cal A}$ is Laplacian operator, consider a positive harmonic function
\[ h_1(x)=log|x| \]
which vanishes on $\partial D_1$. 
For the general operator choose $h_1$ as a Green function of $D_1$ with the pole at $\infty$.
From \cite{LSW} it follows that for arbitrary $\delta>1$ there are positive constants $c_1$ and $c_2$ such that
\begin{equation}\label{Eq1:2}
c_1 log|x| \le h_1(x) \le c_2 log|x|, \quad\text{for}\ x\in D_\delta=\{x\in \rr^2: |x|\geq \delta \}
\end{equation}
Moreover, $h_1$ is H\"{o}lder continuous in $D_1$ and continuously vanishes on $\partial D_1$. Obviously, the Dirichlet problem on $D_1$ has infinitely many solutions 
in a class $O(h_1)$, or otherwise speaking, the logarithmic singularity
at $\infty$ is not removable.

Let $\Omega \subset D_1$ be any open subset and $g \in C(\partial \Omega)$ and
\begin{equation}\label{Eq1:3}
g(x)=O(h_1(x)) \quad\text{as}\ x\in \partial \Omega, \ |x|\to +\infty
\end{equation}
The main goal of this paper is to derive a necessary and sufficient condition for the 
unique solvability of the DP
\begin{equation}\label{Eq1:4}
{\cal A}u = 0 \quad\text{in}\  \ \Omega, \ \ u=g \quad\text{on}\ \ \partial\Omega
\end{equation}
in a class
\begin{equation}\label{Eq1:5}
u(x)=O(h_1(x)) \quad\text{as}\ x\in  \Omega, \ |x|\to +\infty
\end{equation}
without further specification of the behaviour of solution at infinity. 
A typical example of this kind is obtained by taking $g=h_1$ as a boundary function.
It is then clear that $h_1$ is a solution of the DP (\ref{Eq1:4}),(\ref{Eq1:5}), however it is not
clear whether $h_1$ is the only solution. Another interesting example is obtained by prescribing
$g=0$. In this case, $u\equiv 0$ is the only bounded solution, while it is not clear whether it is
the only solution in a class (\ref{Eq1:5}). Otherwise speaking, what is the criteria on $\Omega$ for the
removability of the logarithmic singularity at $\infty$. Equivalently, whether or not zero is an eigenvalue of ${\cal A}$ in $\Omega$ in a class (\ref{Eq1:5}).

An equivalent problem can be posed near the finite boundary point which is the inversion of $\infty$.
Under the inversion $x\to x|x|^{-2}$, $D_1$ is mapped to
\[ D= \{x\in \rr^2: 0<|x|<1 \}\]
If ${\cal A}$ is Laplacian operator, consider the positive harmonic function
\[ h(x)=log\frac{1}{|x|} \]
which vanishes on $\partial D\setminus \{\zeta\}$, where $\zeta$ is the origin. For the general operator ${\cal A}$ choose $h$ as a Green function of $D$ with the pole at $\zeta$.
From \cite{LSW} it follows that for an arbitrary compact subset $E$ of the unit disk $\Sigma= D\cup \{\zeta\}$ there is a constant $C\ge 1$ such that
\begin{equation}\label{Eq1:6}
C^{-1} log\frac{1}{|x|} \le h(x) \le C log\frac{1}{|x|}, \quad\text{for}\ x\in E
\end{equation}
Moreover, $h$ is H\"{o}lder continuous on $D$ and vanishes continuously on $\partial D \setminus \{\zeta\}$. The Dirichlet problem on $D$ has infinitely many solutions 
in the class $O(h)$, or otherwise speaking, the logarithmic singularity
at $\zeta$ is not removable.

Let $\Omega \subset D$ be any open subset and $g \in C(\partial \Omega\setminus \{\zeta\})$ and
\begin{equation}\label{Eq1:7}
g(x)=O(h(x)) \quad\text{as}\ x\in \partial \Omega\setminus \{\zeta\}, \ x\to \zeta
\end{equation}
The main goal is to derive a necessary and sufficient condition for the 
unique solvability of the DP
\begin{equation}\label{Eq1:8}
{\cal A}u = 0 \quad\text{in}\  \ \Omega, \ \ u=g \quad\text{on}\ \ \partial\Omega\setminus \{\zeta\}
\end{equation}
in a class
\begin{equation}\label{Eq1:9}
u(x)=O(h(x)) \quad\text{as}\ x\in  \Omega, \ x\to \zeta
\end{equation}
without further specification of the behaviour of solution at $\zeta$. Equivalently, the goal is to derive necessary and sufficient conditions on $\Omega$ for the removability of the logarithmic singularity at $\zeta$. The equivalence of described problems at $\infty$ and at $\zeta$ follow from Kelvin transformations in the case of a Laplacian operator, and from the generalized Kelvin transformation in the case of general elliptic operator ${\cal A}$. A generalization  
of the well-known Kelvin transformation for the Laplace equation to the equation (\ref{Eq:W:1:1}) was introduced  \cite{SerrinW}. This transformation allows one to transform the uniformly elliptic equation (\ref{Eq:W:1:1}) near $\infty$ to another uniformly elliptic equation with bounded measurable coefficients
near a finite boundary point $\zeta$ which is the inversion of $\infty$. Throughout the rest of the paper, without loss of generality we are going to consider the case of the finite boundary point.
\subsection{Formulation of Problems}\label{E:1:2}
Let $\Omega \subset D$ denote any 
open 
subset, $\pl \Omega$ its topological boundary and $\zeta\in \partial \Omega$.
We consider the differential operator ${\cal A}u$, with 
$a_{ij}=a_{ji}$ being real bounded measurable functions defined
in $\Omega$. 

Throughout this paper, we use the summation convention and assume that
${\cal A}$ is uniformly elliptic. That is, there is a constant 
$\lambda \ge 1$, such that
\begin{equation}\label{Eq:W:1:2}
\lambda^{-1}|\xi|^2 \le 
a_{ij}(x)\xi_i\xi_j\le \lambda |\xi|^2
\end{equation}
for all $x\in \Omega$ and all $\xi \in \rr^N$. We will assume that 
the coefficients $a_{ij}$ are defined and satisfy (\ref{Eq:W:1:2})
for all $x\in D$. This can always be achieved 
by putting $a_{ij}=\delta_{ij}$ outside of $\Omega$. 

Throughout, we use the standard notation of Sobolev spaces \cite{Adams}.
A function $u$ in $H_{loc}^{1,2}(\Omega)$ is a weak solution
of the equation (\ref{Eq:W:1:1}) in $\Omega$ if
\begin{equation}\label{Eq:W:1:3}
\int_\Omega a_{ij}u_{x_i}\phi_{x_j}dx=0
\end{equation}
whenever $\phi \in C_0^{\infty}(\Omega)$. A function $u$ in $H_{loc}^{1,2}(\Omega)$ is a supersolution
of (\ref{Eq:W:1:1}) in $\Omega$ if
\begin{equation}\label{Eq:W:1:4}
\int_\Omega a_{ij}u_{x_i}\phi_{x_j}dx\ge 0
\end{equation}
whenever $\phi \in C_0^{\infty}(\Omega)$ is nonnegative.
A function $u$ is a subsolution of (\ref{Eq:W:1:1}) in $\Omega$ if
$-u$ is a supersolution of (\ref{Eq:W:1:1}). 

The weak solution of (\ref{Eq:W:1:1})
is locally H\"{o}lder continuous \cite{DeGiorgi, Nash, Moser}.
A continuous weak solution of (\ref{Eq:W:1:1}) in $\Omega$ 
is called ${\cal A}$-{\it harmonic} in $\Omega$. 

A function $u$ is called a ${\cal A}$-{\it superharmonic} in $\Omega$ if it satisfies the following conditions:
\begin{description}
\item{\bf(a)} $-\infty<u\le+\infty$, $u < +\infty$ on a dense subset of $\Omega$; 
\item{\bf(b)} $u$ is lower semicontinuous (l.s.c.);
\item{\bf(c)} for each open $U\Subset\Omega$ and each ${\cal A}$-harmonic $h\in C(\ov{U})$,
the inequality $u\ge h$ on $\partial U$ implies $u \ge h$ in $U$.    
\end{description}
In the case when ${\cal A}=\Delta$,
${\cal A}$-harmonic (respectively ${\cal A}$-superharmonic, ${\cal A}$-subharmonic) functions coincide with classical harmonic (respectively superharmonic, subharmonic) functions.

A function $u=v/h$ is called a ${\cal A}_h$-{\it harmonic}, ${\cal A}_h$-{\it superharmonic}, or\\ ${\cal A}_h$-{\it subharmonic} in $\Omega$ if $v$ is ${\cal A}$-harmonic, ${\cal A}$-superharmonic, or ${\cal A}$-subharmonic. In the case when ${\cal A}=\Delta$,
${\cal A}_h$-harmonic (respectively ${\cal A}_h$-superharmonic, ${\cal A}_h$-subharmonic) functions coincide with classical $h$-harmonic (respectively $h$-superharmonic, $h$-subharmonic) functions \cite{Doob}.

We use the notation
${\cal S}(\Omega)$ for a class of all ${\cal A}_h$-superharmonic functions in $\Omega$.
Similarly, $u$ is ${\cal A}_h$-{\it subharmonic} in $\Omega$ if $-u$ is 
${\cal A}_h$-superharmonic in $\Omega$; the class of all ${\cal A}_h$-subharmonic
functions in $\Omega$ is $-{\cal S}(\Omega)$. 

It is well known (\cite{LSW, HKM}) that $u\in H_{loc}^{1,2}(\Omega)$ is 
${\cal A}$-superharmonic if and only if
it is a supersolution with
\begin{equation}\label{Eq:W:1:5}
u(x)=\liminf_{y\to x} u(y)=ess \liminf_{y\to x} \tilde{u}(y) \ \quad\text{for all}\ x \in \Omega. 
\end{equation}
where essential limit inferior means the largest limit inferior in the $L_1$-equivalence
class of functions $\tilde{u}$ such that $\tilde{u}=u$ almost everywhere. Moreover, in the $L_1$-equivalence
class of every supersolution, there is an ${\cal A}$-superharmonic representative
which satisfies (\ref{Eq:W:1:5}).

Given boundary function $f$ on $\partial \Omega$, consider a {\bf $h$-Dirichlet
problem ($h$-DP): find ${\cal A}_h$-harmonic function $u$ in $\Omega$ such that }
\begin{equation}\label{Eq:W:1:6}
u=f \quad\text{on}\ \pl \Omega
\end{equation}
It is easy to see that ${\cal A}_h$-harmonic function $u=\frac{v}{h}$ is a bounded solution of the
$h$-DP if and only if $v$ is a solution of the DP (\ref{Eq1:7})-(\ref{Eq1:9}).

Assuming for a moment
that $f\in C(\partial \Omega)$, the {\it generalized upper (or lower) solution} (in the sense of Perron-Wiener-Brelot)  of the $h$-DP is defined as
\begin{equation}\label{Eq:W:1:7}
^h \bar{H}_f^{\Omega} \equiv \inf \{u\in {\cal S}(\Omega): \liminf_{x\to y, x\in \Omega} u \ge f(y) \ \quad\text{for all}\ y \in \partial\Omega \}
\end{equation}
\begin{equation}\label{Eq:W:1:8}
^h \b{H}_f^{\Omega} \equiv \sup \{u\in {-\cal S}(\Omega): \limsup_{x\to y, x\in \Omega} u \le f(y) \ \quad\text{for all}\ y \in \partial \Omega \}
\end{equation}
The class of functions defined in (\ref{Eq:W:1:7}) (or in (\ref{Eq:W:1:8}) is called {\it upper class} (or {\it lower class}) of the $h$-DP.

The boundary function $f$ is called an {\it $h$-resolutive boundary function} if (\cite{Doob}),  
\[ ^h\bar{H}_f^{\Omega} \equiv \ ^h\b{H}_f^{\Omega} \equiv \ ^hH_f^{\Omega}, \]
and $^hH_f^{\Omega}$ is called a generalized solution of the $h$-DP for $f$.

A boundary of the Greenian open set is called {\it $h$-resolutive} if every finite-valued continuous boundary function is $h$-resolutive, equivalently, if the bounded Borel measurable boundary functions are $h$-resolutive (\cite{Doob}). 

Equivalently, we can define a generalized solution of the DP (\ref{Eq1:7})-(\ref{Eq1:9}):

{\bf Definition 1.1.}\ Let $g$ be a Borel measurable boundary function defined on $\partial \Omega \setminus \{\zeta\}$ that satisfies \eqref{Eq1:7}. Then $g$ is called a resolutive boundary function for the DP  (\ref{Eq1:7})-(\ref{Eq1:9}), if $f=g/h$ (extended to $\zeta$) is $h$-resolutive for the $h$-DP. The function
\begin{equation}\label{newPWB}
H_g^\Omega := h \ ^hH_f^{\Omega}
\end{equation}
is called a generalized solution of the DP (\ref{Eq1:7})-(\ref{Eq1:9}).

For a given boundary Borel subset $A\subset \partial \Omega$, denote the indicator function
of $A$ as $1_A$. If $1_A$ is a resolutive boundary function we define the {\it $h$-harmonic measure} of $A$ as (\cite{Doob}):

\[ \mu_{\Omega}^h(\cdot,A)= \ ^hH_{1_A}^\Omega(\cdot). \]

It is said that $A$ is an ${\cal A}_h$-harmonic measure null set if $\mu_{\Omega}(\cdot,A)$
vanishes identically in $\Omega$. If this is not the case, $A$ is a set of
positive ${\cal A}_h$-harmonic measure. Since $^hH_{1_{A}}^\Omega$ is ${\cal A}_h$-harmonic
in $\Omega$, from the strong maximum principle it follows that if $\mu_{\Omega}^h(x,A)>0$,
then it is positive in the whole connected component of $\Omega$ which contains $x$. 

It is easy to see that the boundary of $D$ is $h$-resolutive and the only possible solutions of the $h$-DP in $D$ are constants. Precisely, the unique solution of the $h$-Dirichlet Problem is identical with the constant $f(\zeta)$. Indeed, for arbitrary $\epsilon >0$, the function 
\[ u(\cdot)=f(\zeta)+\frac{\epsilon}{h(\cdot)} \ \ \Big (\quad\text{or}\ v(\cdot)=f(\zeta)-\frac{\epsilon}{h(\cdot)} \Big ) \]
is in the upper class (or lower class) for $h$-DP in $D$ for $f$. Hence,
\[ f(\zeta)-\frac{\epsilon}{h(\cdot)}\leq \ ^h \b{H}_f^{D}(\cdot) \leq \ ^h \bar{H}_f^{D}(\cdot) \leq f(\zeta)+\frac{\epsilon}{h(\cdot)} \]
Since, $\epsilon>0$ is arbitrary, the assertion follows.
In particular, we have
\begin{equation}\label{Eq:W:1:8a}
\mu_{D}^h(\cdot,\partial D -\{\zeta\})\equiv 0, \ \ \mu_{D}^h(\cdot,\{\zeta\})\equiv 1.
\end{equation}
Since $\partial D$ is $h$-resolutive, from \cite{Doob} (see Example(b) and its continuation, pages 116-117) it follows that $\partial \Omega$ is $h$-resolutive for $h$-DP in $\Omega$. It is not difficult to see that 
any Borel subset $A\subset \partial\Omega \cap (\partial D -\{\zeta\})$ is an $h$-harmonic measure null set in $\Omega$.
Indeed, for arbitrary $\epsilon>0$, $\epsilon h$ (or $-\epsilon h$) is in the upper class (or lower class) of the $h$-DP in $\Omega$ with boundary function $1_A$. This fact also follows from the following formula proved in \cite{Doob} (formula (8.3), p.117): if $A$ is a Borel subset of $\partial D$, then
\begin{equation}\label{Eq:W:1:8b}
\mu_{D}^h(\xi,A)=\mu_{\Omega}^h(\xi,A\cap \partial\Omega)+\int\limits_{D\cap \partial\Omega} \mu_{D}^h(\eta,A)\mu_{\Omega}^h(\xi,d\eta), \ \ \xi\in\Omega
\end{equation}
By choosing $A=\partial D -\{\zeta\}$, from (\ref{Eq:W:1:8a}),(\ref{Eq:W:1:8b}) it follows that 
\begin{equation}\label{Eq:W:1:8c}
\mu_{\Omega}^h\Big(\xi,\partial\Omega\cap(\partial D -\{\zeta\})\Big)= \ ^hH_{1_{\partial\Omega\cap(\partial D -\{\zeta\})}}^\Omega(\xi) \equiv 0, \ \ \xi\in\Omega
\end{equation}
By choosing $A=\{\zeta\}$, from (\ref{Eq:W:1:8a}),(\ref{Eq:W:1:8b}) it follows that 
\begin{equation}\label{Eq:W:1:8d}
1=\mu_{\Omega}^h(\xi,\{\zeta\})+\mu_{\Omega}^h(\xi,D\cap\partial\Omega), \ \ \xi\in\Omega
\end{equation}
or
\begin{equation}\label{Eq:W:1:8e}
1= \ ^hH_{1_{\{\zeta\}}}^\Omega(\xi)+ \ ^hH_{1_{D\cap\partial\Omega}}^\Omega(\xi), \ \ \xi\in\Omega
\end{equation}
All the proofs referenced for the Laplacian operator remain the same for the general elliptic operator (\ref{Eq:W:1:1}). Hence, the indicator function of any Borel measurable boundary subset of $\partial\Omega$, and equivalently any bounded Borel measurable boundary function is resolutive for the $h$-DP in $\Omega$. 
The following formula is true for the solution  $^h H_f^{\Omega}$ of the $h$-Dirichlet Problem:
\begin{equation}\label{Eq:W:1:8f}
^h H_f^{\Omega}(\xi)=\int\limits_{\Sigma\cap\partial\Omega}f(\eta)\mu_\Omega^h(\xi,d\eta), \ \ \xi\in\Omega
\end{equation}
where $\Sigma =D\cup \{\zeta\}$ be a unit disk. In fact, the class of $h$-resolutive boundary functions is larger than the class of bounded Borel measurable functions. Precisely, the boundary function $f$ is $h$-resolutive if and only if $f 1_\Sigma \in L^1(\mu_\Omega^h)$, where the latter denotes the class of $\mu_\Omega^h$-integrable boundary functions on $\Sigma\cap\partial\Omega$. By using slight modification of the notation, (\ref{Eq:W:1:8f}) will be written in the following equivalent form:
\begin{equation}\label{Eq:W:1:8h}
^h H_f^{\Omega}(\xi)=\mu_\Omega^h(\xi,f 1_\Sigma), \ \ \xi\in\Omega
\end{equation}
The generalized solution is unique by construction, and it coincides with the classical, or Sobolev space solution whenever the latter exists. Also,
note that the construction of the generalized solution of the $h$-DP 
is accomplished by prescribing the behavior of
the solution at $\zeta$. 

Equivalently, according to the formula \eqref{newPWB} the following formula is true for the unique solution of the DP (\ref{Eq1:7})-(\ref{Eq1:9}):
\begin{equation}\label{Eq:W:1:8h}
H_g^{\Omega}(\xi)=h(\xi)\int\limits_{\Sigma\cap\partial\Omega}\frac{g(\eta)}{h(\eta)}\mu_\Omega^h(\xi,d\eta), \ \ \xi\in\Omega
\end{equation}
Again note that the construction of the unique solution $H_g^{\Omega}$ of the DP (\ref{Eq1:7})-(\ref{Eq1:9}) is accomplished by prescribing the behaviour 
of the ratio $H_g^{\Omega} /h$ at $\zeta$.
The elegant theory, 
while identifying a class of unique solvability, 
leaves the following questions open:
\begin{itemize}
\item Would a unique solution of the $h$-DP still exist if its 
limit at $\zeta$ were not specified? That is, could it be 
that the solutions would pick up the ``boundary value" $f(\zeta)$
without being required? Equivalently, would unique solution of the DP (\ref{Eq1:7})-(\ref{Eq1:9}) still exist if the limit of the ratio of solution to $h$ at $\zeta$ is not prescribed?  In particular, is the logarithmic singularity at $\zeta$
removable?  
\item What if the boundary datum $f$ (or $g/h$) on $\pl \Omega$, while being continuous at $\partial\Omega\setminus\{\zeta\}$, 
does not have a limit at $\zeta$, for example, it 
exhibits bounded oscillations. Is the $h$-Dirichlet Problem (or DP (\ref{Eq1:7})-(\ref{Eq1:9})) uniquely solvable?
\end{itemize}
The example given above demonstrates that if $\Omega= D$, the answer is negative and 
arbitrary constant $C$ is a solution of the $h$-Dirichlet problem, $Ch$ is a solution of the DP (\ref{Eq1:7})-(\ref{Eq1:9}) and that the logarithmic singularity at $\zeta$ is not removable. The positive answer to these fundamental questions is possible if 
$\Omega$ is not too sparse, or equivalently $\Omega^c \cap D$ is
not too thin near $\zeta$. 
The principal purpose of this paper is to prove  
the Wiener test for the thinness of $\Omega^c \cap D$ near $\zeta$ 
which guarantees the uniqueness of the solution of the
$h$-Dirichlet Problem (or DP (\ref{Eq1:7})-(\ref{Eq1:9})) without specification of the boundary function (or ratio of the boundary function to $h$) at $\zeta$.

Furthermore, we assume that $f=g/h: \partial\Omega\setminus\{\zeta\} \rightarrow \rr$ 
is a {\it bounded Borel measurable} function. 
By fixing an arbitrary finite real number $\bar{f}$, extend a function
$f$ as $f(\zeta)=\bar{f}$. Obviously, the extended function
is a bounded Borel measurable on $\partial\Omega$ and
there exists a unique solution $^hH_f^{\Omega}$ of the $h$-DP, and the unique solution of the DP (\ref{Eq1:7})-(\ref{Eq1:9})
is given by \eqref{newPWB}. The major question now becomes:\\

{\it {\bf Problem 1:} How many bounded solutions do we actually have, or does the constructed
solution depend on} $\bar{f}$ ? \\


Since the indicator function $1_{\{\zeta\}}$ is a resolutive boundary function,
the ${\cal A}_h$-harmonic measure of $\{\zeta\}$ is well defined:

\[ \mu_{\Omega}^h(\cdot,\{\zeta\})= \ ^hH_{1_{\{\zeta\}}}^\Omega(\cdot). \]

Next problem is the measure-theoretical counterpart
of the Problem 1:\\

{\it {\bf Problem 2:} Given $\Omega$, is the ${\cal A}_h$-harmonic measure 
of $\{\zeta\}$ null or positive} ?\\


In fact, both major problems are equivalent, and the next definition expresses the 
connection between them.\\  

{\bf Definition 1.2.}\ $\zeta$ is said to be $log$-{\it regular} for
$\Omega$ if it is an ${\cal A}_h$-harmonic measure null set. Conversely, $\zeta$ is $log$-{\it irregular} if it has a positive ${\cal A}_h$-harmonic measure.\\

The notion of the regularity of $\zeta$ is, in particular, related to 
the notion of continuity of the solution at $\zeta$. 

{\it {\bf Problem 3:} Given $\Omega$, such that $\zeta$ is the limit point of $\partial\Omega$,
whether or not}

\[\liminf_{x\to \zeta, x\in \partial\Omega} f \le \liminf_{x\to \zeta, x\in \Omega}  \ ^hH_f^{\Omega} \le \limsup_{x\to \zeta, x\in \Omega}  \ ^hH_f^{\Omega} \le \limsup_{x\to \zeta, x\in \partial\Omega} f \]
\begin{equation}\label{Eq:W:1:10}
 \quad\text{for all bounded}\ f\in C(\partial \Omega\setminus\{\zeta\}).
\end{equation}

Note that if $f$ has a limit at $\zeta$, (\ref{Eq:W:1:10}) simply means that 
the solution $^hH_f^{\Omega}$ is continuous at $\zeta$.

The equivalent problem in the context of the DP (\ref{Eq1:7})-(\ref{Eq1:9}) is the following:

{\it {\bf Problem 3$^\prime$:} Given $\Omega$, such that $\zeta$ is the limit point of $\partial\Omega$,
whether or not}

\[\liminf_{x\to \zeta, x\in \partial\Omega} \frac{g}{h} \le \liminf_{x\to \zeta, x\in \Omega}  \ \frac{H_g^{\Omega}}{h} \le \limsup_{x\to \zeta, x\in \Omega}  \ \frac{H_g^{\Omega}}{h} \le \limsup_{x\to \zeta, x\in \partial\Omega} \frac{g}{h} \]
\begin{equation}\label{Eq:W:1:10aa}
 \quad\text{for all}\ g \quad\text{such that}\ \frac{g}{h}\in C(\partial \Omega\setminus\{\zeta\}) \quad\text{and bounded}.
\end{equation}

In particular, if  $g/h$ has a limit at $\zeta$, (\ref{Eq:W:1:10aa}) means that 
the limit of the ratio $H_g^{\Omega}/h$ at $\zeta$ exists and equal to the limit of $g/h$.\\

The notion of $log$-{\it regularity} of $\zeta$ introduced in Definition 1.2 fits naturally in the framework of minimal-fine topology. Recall that ${\cal A}$-fine topology
is the coarsest topology of $\rr^2$ which makes every superharmonic function continuous
\cite{Doob, HKM}.  ${\cal A}$-fine topology is finer than the Euclidean topology. It is well-known that there is an elegant connection
between a problem of finding the structure of the neigborhood base in
${\cal A}$-fine topology and the problem on the regularity of finite boundary points.
Namely, given open set $\Omega \subset D$, its boundary point $x_0\in D$
is irregular if and only if $\Omega$ is a deleted neigborhood of $x_0$
in ${\cal A}$-fine topology \cite{Doob, Landkov, HKM}. Equivalently, $\Omega^c$ is called thin at 
$x_0$. Hence, for any open subset $\Omega \subset D$ with $\zeta \in \partial\Omega$, thinness 
of $\Omega^c$ at $\zeta$ in ${\cal A}$-fine topology is equivalent to the classical regularity 
of the boundary point $\zeta$, which is characterized through the celebrated Wiener test for the boundary regularity of harmonic functions (\cite{Wiener1,Wiener2}) and it is independent of the elliptic operator ${\cal A}$ \cite{LSW}. 
One can also introduce ${\cal A}$-minimal-fine topology at $\zeta$ as an extension of the ${\cal A}$-fine topology of $D$ to isolated boundary point $\zeta$. Positive ${\cal A}$-harmonic function $h$ is minimal ${\cal A}$-harmonic in $D$, in the sense that it dominates there no positive ${\cal A}$-harmonic function except for its own constant submultiples. $\zeta$ is a minimal Martin boundary point of $D$ and $h$ is an associated minimal Martin boundary function
with pole at $\zeta$. 
${\cal A}$-minimal-fine topology has as relative topology on $D$ the fine topology. \\

{\bf Definition 1.3.}\ Subset $E\subset D$ is ${\cal A}$-minimally thin at $\zeta$ if
\begin{equation}\label{Eq:W:1:10a}
\lim_{B\downarrow \zeta} \hat{R}^{E\cap B}_h(\xi)=0, \ \xi\in D
\end{equation}
where $\lim$ as $B\downarrow \zeta$ means the Euclidean neigborhood of $\zeta$ shrinks to $\zeta$, that is, as its diameter shrinks to 0; $R^{E}_h$ becomes a reduction of $h$ on E
\[ R^{E}_h(x)=\inf \{v(x): v \quad\text{is ${\cal A}$-superharmonic in} \ \Omega, v\ge 0, v\ge h \quad\text{on} \ E \} \]
and $\hat{R}^{E}_h$ is its lower-semicontinuous regularization, or smoothed reduction:
\[ \hat{R}^{E}_h(x)=\lim_{\epsilon \to 0} \inf_{y\in B_\epsilon(x)}R^{E}_h(y) \]

Note that the sequence  $\hat{R}^{E\cap B}_h(\xi)$ is monotonically decreasing as $B\downarrow \zeta$, and it is satisfactory if \eqref{Eq:W:1:10a}
is satisfied for one particular choice of $B\downarrow \zeta$.

The following is the equivalent definition of ${\cal A}$-minimally thinness at $\zeta$ (see Lemma \ref{minimalthinness1}).\\

{\bf Definition 1.3$^\prime$.} \ Subset $E\subset D$ is ${\cal A}$-minimally thin at $\zeta$ if
\begin{equation}\label{minimalthinnessdefinition}
\hat{R}^{E\cap B}_h(\xi)\not\equiv h(\xi), \ \xi\in D
\end{equation}
for some (and equivalently for all) Euclidean neigborhood(s) $B$ of $\zeta$.


A point $\zeta$ is a minimal-fine limit point of a set $E\subset D$ if $E$ is not ${\cal A}$-minimal thin at $\zeta$. We write this fact as $\zeta \in E^{mf}$, where $E^{mf}$ is the set of minimal-fine limit points of $E$. If $E$ is ${\cal A}$-minimally thin at $\zeta$, then equivalently, $E^c\cap D$ is a deleted ${\cal A}$-minimal-fine neighborhood of $\zeta$. In fact, ${\cal A}$-fine topology and ${\cal A}$-minimal-fine topology are equivalent at the isolated boundary point $\zeta$. 

We can now formulate the topological counterpart of problems 1, 2 and 3:\\

{\it {\bf Problem 4:} Is the given open set $\Omega\subset D$ a deleted neigborhood of $\zeta$
in ${\cal A}$-minimal-fine topology? Equivalently, is $\Omega^c$ ${\cal A}$-minimally thin at $\zeta$?
Or conversely, whether or not $\zeta \in (\Omega^c\cap D)^{mf}$.} \\

The principal result of this paper expresses the solutions to equivalent
Problems 1-4 in terms of the Wiener test for the $log$-{\it regularity} of $\zeta$ (or $\infty$).
\subsection{The Main Result}\label{mainresult}

Key for solution of the equivalent Problems 1-4 is the notion of $h$-Capacity. We introduce and analyze this fundamental notion in Section \ref{E:1:3}. To formulate the main result, let us introduce the notion of $h$-Capacity of compact subsets of $D$.\\  

{\bf Definition 1.4.}\ The $h$-Capacity of the compact $K \subset D$ with respect to the elliptic operator ${\cal A}$ is the number
\[ C^{\cal A}_h(K)= \inf_{\phi \in  {\cal G}} {\cal L}_{\cal A}(\phi) \]
where
\[ {\cal L}_{\cal A}(\phi) = \int_D a_{ij}\phi_{x_i}\phi_{x_j} dx  \]
\[ {\cal G} =  \Big \{ \phi \in H_0^{1,2}(D): \phi \geq h \quad\text{on} \ K  \quad\text{in the sense of} \ H_0^{1,2}(D) \Big \} \]
where $h(\cdot)=G_{{\cal A}}(\cdot, \zeta)$, and $G_{{\cal A}}$ is the Green function of the operator ${\cal A}$ on unit disk $\Sigma$ with center at origin $\zeta$.\\

 The following is the equivalent definition of $h$-Capacity:\\

{\bf Definition 1.5.}\ If $K\subset D$ is compact, the $h$-Capacity of $K$ is
\[ C^{{\cal A}}_h(K)\equiv sup\{\mu(K): \mu \in {\cal M}(K), \ U_{{\cal A}}^\mu \le h \quad\text{in}\ D \}, \]
where ${\cal M}(K)$ denotes the collection of all nonnegative Radon measures with
support $S(\mu)\subseteq K$; \ $U_{{\cal A}}^\mu$ is the Greenian potential of the form
\[ U_{{\cal A}}^\mu(x) = \int\limits_{K} \frac{G_{{\cal A}}(x,y)}{h(y)}\mu(dy).  \]

We omit index ${\cal A}$ and simply write  $C_h$, $U^\mu$ and $G$  if ${\cal A}=\Delta$ is Laplacian.
Note that $G$ is the Green's function of the Laplacian operator in a unit disk:
\begin{equation}\label{Greenfunction}
G(\xi,\eta)=\log \frac{|\xi^\prime - \eta||\xi|}{|\xi - \eta|} 
\end{equation}
where $\xi^\prime$ is the image of $\xi$ under inversion in the unit circle.\\

In Section \ref{E:1:3} we demonstrate that there is a unique measure $\gamma \in {\cal M}(K)$ such that $C^{{\cal A}}_h(K)=\gamma(K)$
and the potential $U_{{\cal A}}^\gamma$ coincides with the smoothed reduction of $h$ on $K$:
\[U_{{\cal A}}^{\gamma}(x)= \hat{R}^{K}_h(x)=\int\limits_{K} \frac{G_{{\cal A}}(x,y)}{h(y)}\gamma(dy) \]
Measure $\gamma$ is called the $h$-equilibrium measure of $K$ and the Green potential $U_{{\cal A}}^\gamma$ is called the $h$-equilibrium potential.

If ${\cal A} = \Delta$ and $h\equiv 1$, then $C_h(K)$ coincides with the Greenian capacity, the related $\gamma$ becomes the Greenian equilibrium
measure, and $U^\gamma$ becomes the Greenian equilibrium potential, which is the smoothed reduction of 1. We will denote Greenian capacity as $C_g(K)$:
\[ h\equiv 1 \ \Rightarrow \ C_g(K)=C_h(K), \  U^{\gamma}(x)= \hat{R}^{K}_1(x)=\int\limits_{K} G(x,y)\gamma(dy)  \]

Let $a>1$ is fixed and 
\[ E_n=\Omega^c\cap\{\xi: a^{n}\le h(\xi)\le a^{n+1} \}  \]
Our main theorem reads:
\begin{theorem}\label{main theorem}
Given open set $\Omega \subset D$, the following conditions are equivalent:
\begin{description}
\item{(1)\ } $\zeta$ is $log$-regular (or $log$-irregular).
\item{(2)\ }The DP (\ref{Eq1:7})-(\ref{Eq1:9}) has a unique (or infinitely many) solution(s) in a class $O(\log |\cdot|)$.
\item{(3)\ }Boundary regularity conditions (\ref{Eq:W:1:10}),(\ref{Eq:W:1:10aa}) are satisfied (respectively aren't satisfied).
\item{(4)\ }$\Omega^c$ is not ${\cal A}$-minimally thin at $\zeta$ (or it is)
\item{(5)\ }Both of the following series diverge (or converge)
\begin{equation}\label{Eq:W:1:12}
\sum_{n} a^{-n} C_h^{{\cal A}}(E_n), \ \ \  \sum_{n} a^{-n} \hat{R}^{E_n}_h(\zeta)
\end{equation}
\item{(6)\ }The following integral diverges (or converges)
\begin{equation}\label{Eq:W:1:13}
\int_1^{+\infty} \frac{c(\rho)}{\rho^2}d\rho
\end{equation}
where
\[ c(\rho)=C_h^{{\cal A}}(E^\rho), \ E^\rho = \Omega^c \cap \{\xi: 1 \leq h(\xi) \leq \rho \}. \]
\item{(7)\ }Both series \eqref{Eq:W:1:12} and the integral \eqref{Eq:W:1:13} diverge (or converge) when ${\cal A}$ is a Laplace operator.
\item{(8)\ }Both of the following series diverge (or converge) 
\begin{equation}\label{Eq:W:1:14}
\sum_{n}a^n C_g(E_n), \ \ \ \sum_{n} \hat{R}^{E_n}_1(\zeta)
\end{equation}
\end{description}
\end{theorem}

The claim $(8)$ of Theorem \ref{main theorem} is a classical Wiener test for the boundary regularity of harminic functions \cite{Wiener1, Wiener2}. 
The same criterion hold for the boundary regularity of solutions of the second order uniformly elliptic PDE \eqref{Eq:W:1:1} \cite{LSW}. Hence, Theorem \ref{main theorem} establishes that the criterion for the removability of the logarithmic singularity and uniqueness in $O(\log)$ class for the second order uniformly elliptic PDEs with bounded measuarable coefficients coincides with the classical Wiener criterion.

Theorem \ref{main theorem} implies the following counterexample. Let 
\[ \Omega = D \setminus \{(x,0)\in R^2: 0<x<1\} \]
be a unit disk with deleted radius. Then the Wiener series is divergent and the logarithmic singularity at $\zeta$ is removable for any uniformly elliptic operator ${\cal A}$. 

Let $a>1,N>1$ are fixed, and consider
\[ \Omega=D\setminus \bigcup_{n=N}^{+\infty} \{(x,0)\in R^2: e^{-a^{n+1}} \leq x \leq  e^{-a^{n+1}}+\delta_n \} \]
where
\[ \delta_n=e^{-a^{n}n \log n \log_2 n \cdots \log_k^{1+\epsilon} n}, \ k>1\quad\text{is any integer} \] 
Then the Wiener series diverges or converges, that is to say the logarithmic singularity either is or isn't removable, depending on whether $\epsilon \leq 0$ or $\epsilon > 0$. 

{\bf Probabilistic Counterpart:} From the probabilistic standpoint, the Wiener test 
of Theorem ~\ref{main theorem} presents an asymptotic probability laws for conditional 
Markov processes. Although we are not going to present probabilistic proofs in this paper, let us formulate asymptotic laws in the case of the Laplacian operator. Let us consider $h$-Brownian motion on $D$, where $h(\xi)=-\log|\xi|$ and $\zeta$ is the notation for the origin. It is proven in \cite{Doob} (Theorem 2.X.4) that almost every $h$-Brownian path $w_\xi^h(t)$ from $\xi\in D$ has a finite lifetime $S_\xi^h$ and tends to $\zeta$ at its lifetime. 
Let ${\bf B}$ be an event that $\{t: w_\xi^h(t)\in \Omega^c\}$ clusters to $S_\xi^h$. Then the {\bf probabilistic counterpart of Theorem \ref{main theorem}} is the following asymptotic law of behaviour of $h$-Brownian path $w_\xi^h(t)$ as $t\to S_\xi^h-0$:
\[ P({\bf B})=0 \quad\text{or}\ 1 \quad\text{according as}\  \sum_{n} a^{-n} C_h(E_n) < \quad\text{or} =+\infty \]
In \cite{Doob} (Section 2.X.5) it is demonstrated that $h$-Brownian motion can start from $\zeta$ - infinity of $h$. In fact, it is proved with the standard procedure of Kolmogorov that $h$-Brownian motion trajectory $w_\zeta^h(t)$ in $D$ from $\zeta$ is a process whose paths have finite lifetimes and tend back to $\zeta$ at these lifetimes. Let
\[ m_{\Omega^c}=\inf\{t: t>0, w_\zeta^h(t) \in \Omega^c \} \]
be the Markovian exit time from $\Omega$. Then the {\bf probabilistic counterpart of Theorem \ref{main theorem}} is the following asymptotic law of behaviour of $h$-Brownian path $w_\zeta^h(t)$ as $t\to 0+$:
\[ P(m_{\Omega^c}=0)=0 \quad\text{or}\ 1 \quad\text{according as}\  \sum_{n} a^{-n} C_h(E_n) < \quad\text{or} =+\infty \]
Similar probability laws can be formulated in the case of general elliptic operator ${\cal A}$, and associated continuous time, time-homogeneous Markov
process with infinitesimal Dynkin generator being a differential operator $-{\cal A}$\cite{Dynkin}. One can introduce conditional $G_{\cal A}(\cdot,\zeta)$-Markov process similar to $h$-Brownian motion. The probabilistic counterpart of the Theorem \ref{main theorem} will be similar asymptotic probability laws characterizing the behaviour of the process as approaching $\zeta$ at the end of lifetime, and local asymptotics of the process starting at $\zeta$, as they are formulated for $h$-Brownian motion, just by replacing $h$-Brownian path with $G_{\cal A}(\cdot,\zeta)$-Markov path.

\section{$h$-Capacity}\label{E:1:3}
In this section we introduce an important, new notion of $h$-{\it capacity}. Let $\Sigma$ be an open unit disk.
For any Radon measure $\mu$ with compact support in $\Sigma$, consider a problem
\begin{equation}\label{Eq:W:2:1}
{\cal A}u=\mu \quad\text{in}\ \Sigma,\ \  u=0 \quad\text{on}\ \pl\Sigma.
\end{equation}

$u\in L^1(\Sigma)$ is a weak solution of (\ref{Eq:W:2:1}) if
\[ \int_\Sigma u {\cal A}\Phi dx = \int_\Sigma \Phi d\mu  \]
for every $\Phi\in H_0^{1,2}(\Sigma)\cap C(\ov{\Sigma})$ such that ${\cal A}\Phi \in C(\ov{\Sigma})$.

The Green's function $G_{{\cal A}}(x,y)$ of the operator ${\cal A}$ on $\Sigma$ is defined as the weak solution 
of the problem (\ref{Eq:W:2:1}) with $\mu=\delta_y$, where $\delta_y$ is the Dirac measure of $y$.
Recall that if ${\cal A}=\Delta$, we omit index ${\cal A}$, and $G$ is the Green's function of the Laplacian in a unit disk given in (\ref{Greenfunction}).

In the next lemma we bring together some well-known facts from \cite{LSW}:

\begin{lemma}\label{2}
\begin{description}
\item{(1)\ } There exists a unique weak solution $u$ of (\ref{Eq:W:2:1}) such that $u\in H_0^{1,p}(\Sigma)$ for every $p<2$.
\item{(2)\ }  Weak solution $u\in H_0^{1,2}(\Sigma)$, if and only if $\mu \in H^{-1,2}(\Sigma)$. 
\item{(3)\ } $G_{{\cal A}}(\cdot,y)\ge 0$ is ${\cal A}$-harmonic and H\"{o}lder continuous
in $\Sigma-y$, $\lim_{x\to y} G_{{\cal A}}(x,y)=+\infty$. 
\item{(4)\ } For every Radon measure with compact support in $\Sigma$, the integral
\begin{equation}\label{Eq:W:2:2}
u(x)=\int_\Sigma G_{{\cal A}}(x,y)\mu(dy).
\end{equation}
exists and finite a.e., and is a weak solution of (\ref{Eq:W:2:1}). 
\item{(5)\ } Let $G_{{\cal A}}$ and $\ov{G}_{{\cal A}}$  be the Green functions for any uniformly elliptic operators ${\cal A}$ and $\ov{{\cal A}}$ with the  ellipticity constant $\lambda$ on a sphere $\Sigma$. Then, for any compact subset $E$ of $\Sigma$, there exists a constant $C\ge 1$ depending on $E, \Sigma$ and $\lambda$ such that
\begin{equation}\label{Eq:W:2:3}
C^{-1}G_{\ov{{\cal A}}}(x,y) \le G_{{\cal A}}(x,y) \le C G_{\ov{{\cal A}}}(x,y), \quad\text{for all}\ x,y\in E.
\end{equation}
\end{description}
\end{lemma} 
If $\ov{\cal A}$ is taken to be the Laplace operator, then from (\ref{Eq:W:2:3}) it follows that for any compact subset $E$ of $\Sigma$,
\begin{equation}\label{Eq:W:2:3a}
C^{-1}log\frac{|\xi^\prime - \eta||\xi|}{|\xi-\eta|} \le  G_{{\cal A}}(\xi,\eta) \le C log\frac{|\xi^\prime - \eta||\xi|}{|\xi-\eta|} \ \quad\text{for all}\ \xi, \eta \in E
\end{equation}
where $\xi^\prime$ is the image of $\xi$ under inversion in unit circle; $C$ depends on $E,\Sigma, \lambda$. Note that $G_{{\cal A}}(\cdot,\zeta)=h(\cdot)$, and in particular, (\ref{Eq1:6}) follows from (\ref{Eq:W:2:3a}) by choosing $\eta=\zeta$. Observe also that there is a constant $C$ depending on $E$ such that
\begin{equation}\label{Eq:W:2:3b}
C^{-1}log\frac{2}{|\xi-\eta|} \le log\frac{|\xi^\prime - \eta||\xi|}{|\xi-\eta|}  \le Clog\frac{2}{|\xi-\eta|} \ \quad\text{for all}\ \xi, \eta \in E
\end{equation}
From (\ref{Eq:W:2:3a}),(\ref{Eq:W:2:3b}) it follows that for any compact subset $E$ of $\Sigma$
\begin{equation}\label{Eq:W:2:3c}
C^{-1}log\frac{2}{|\xi-\eta|} \le G_{{\cal A}}(\xi,\eta)  \le Clog\frac{2}{|\xi-\eta|} \ \quad\text{for all}\ \xi, \eta \in E
\end{equation}
with some constant $C\ge 0$ depending on $E,\Sigma,\lambda$.

Let $K\subset D$ be a compact set and ${\cal M}_K$ be a a set of all nonnegative Radon measures with support on $K$. Introduce the energy integral on ${\cal M}_K$
\[ \Vert \mu \Vert ^2=I_h(\mu)=\iint\limits_{K\times K} G_{\cal A}^h(\xi,\eta)\mu(d\xi)\mu(d\eta), \]
where 
\[ G_{\cal A}^h(\xi,\eta)=\frac{G_{\cal A}(\xi,\eta)}{h(\xi)h(\eta)} \]
Since the symmetric kernel $G_{\cal A}^h$ is nonnegative and lower semicontinuous
in $D\times D$, it follows that $I(\mu)$ is a lower semicontinuous functional in a vague topology (or $*$-weak topology)
of ${\cal M}_K$  \cite{Landkov}. Let
\[ \overline{{\cal M}}_K=\{\mu\in {\cal M}: \mu(K)=1 \}. \]
Since $\overline{{\cal M}}_K$ is convex and vaguely compact, and
\begin{equation}\label{Eq:W:1:11}
 0<h_1=\min\limits_{K}h \le h(x) \le \max\limits_{K}h = h_2 <+\infty
\end{equation}
classical theory (\cite{Landkov}) of equilibrium measure and capacity produced by the Greenian kernel $G(x,y)$ ($h\equiv 1$) 
can be easily extended to the case of general $h$-Greenian kernel $G_h$.  From (\ref{Eq:W:1:11})  it follows that
\[ G_{\cal A}^h(x,y)\ge C_1 log\frac{2}{d}, \ x,y\in K, \]
where $d<2$ is the diameter of $K$. Hence,
\begin{equation}\label{infpositive}
0<W_h(K):=\inf\limits_{\overline{{\cal M}}_K}I_h(\mu) \le +\infty. 
\end{equation}
Minimizing sequence $\{\mu_n\}\subset {\overline{{\cal M}}_K}$ vaguely converges to minimizing measure $\lambda_K\in {\overline{{\cal M}}_K}$, i.e.
\[I_h(\lambda_K)=W_h(K). \]
One can tell more by introducing subclass of measures with finite energy $I_h$:
\[ {\cal E}_K=\{\mu\in {\cal M}: I_h(\mu)<+\infty\}, \ \overline{{\cal E}_K}={\cal E}_K \cap  \overline{{\cal M}}_K \]
In fact, ${\cal E}_K$ has a Hilbert space structure with inner product
\[ (\mu,\nu)=I_h[\mu,\nu]:=\iint\limits_{K\times K} G_{\cal A}^h(\xi,\eta)\mu(d\xi)\nu(d\eta). \]
By classical argument (\cite{Landkov}), one can deduce from (\ref{Eq:W:1:11}) that ${\cal E}_K$ is complete, and the 
variational problem on $\overline{{\cal E}_K}$ has a unique solution $\lambda_K$ which is a strong limit of the minimizing sequence.

{\bf Definition 1.6.}\ The $h$-Capacity of the compact set $K \subset D$ is the number
\[ C^{\cal A}_h(K)=W_h^{-1}(K) \]
We simply write $C_h(K)$ if ${\cal A}=\Delta$.

From (\ref{infpositive}) it follows that $C^{\cal A}_h(K)$ is finite and nonnegative. It coincides with classical Green capacity if $h\equiv 1$. Hence, from (\ref{Eq:W:1:11}) it follows that the compact set $K\subset D$ is of zero Greenian (or Newtonian) capacity if and only if $C^{\cal A}_h(K)=0$. Recall that the set $E$ is said to have inner capacity zero, if any of its compact subsets are of zero capacity and the term ``approximately everywhere" means ``everywhere, except for a set of inner capacity zero". Thus, the concept ``approximately everywhere" (or ``a.e.") for subsets of any compact $K\subset D$ is independent of $h$.  

Consider the following Greenian potential associated with minimizing measure $\lambda=\lambda_K$:
\[ U_{\cal A}^{\lambda}(x)=\int\limits_{K} \frac{G_{\cal A}(x,y)}{h(y)}\lambda(dy)  \]
It has the following properties:
\begin{description}
\item{(1)\ } $ U_{\cal A}^{\lambda}(x)=W_h(K)h(x)=\Vert \lambda \Vert^2 h(x)$ a.e. on the support $S(\lambda)\subset K$ of $\lambda$.
\item{(2)\ } $ U_{\cal A}^{\lambda}(x)\le \Vert \lambda \Vert^2 h(x)$ on  $S(\lambda)$, and accordingly by maximum principle on $D$.
\item{(3)\ } Let $\mu \in \overline{{\cal E}_K}$ and $ U_{\cal A}^{\mu}(x)=W_h(K)h(x)$ a.e. on the support $S(\mu)$. Then $\mu=\lambda$ is a minimizing measure and $\Vert \mu \Vert^2=W_h(K)$.
\item{(4)\ } Let $\mu \in \overline{{\cal E}_K}$ and $ U_{\cal A}^{\mu}(x)=Ah(x)$ a.e. on the support $S(\mu)$, with some constant $A$. Then $\mu=\lambda$ is a minimizing measure and $A=\Vert \mu \Vert^2=W_h(K)$.
\end{description} 
Proof of the properties (1)-(4) coincides with the proof given in (\cite{Landkov}) for the case $h\equiv 1$.\\

{\bf Definition 1.7.}\ Measure
\[ \gamma=W_h^{-1}(K)\lambda=\frac{\lambda}{\Vert \lambda \Vert^2}=C^{\cal A}_h(K)\lambda \]
is called the $h$-equilibrium measure of the compact $K$ and the Green potential $U_{\cal A}^{\gamma}$ is called the $h$-equilibrium potential.\\

From the properties (1),(2) of $U_{\cal A}^\gamma$ it follows that the $h$-equilibrium potential satisfies the properties:
\begin{description}
\item{(1)\ } $ U_{\cal A}^{\gamma}(x)=h(x)$ a.e. on $S(\gamma)\subset K$.
\item{(2)\ } $ U_{\cal A}^{\gamma}(x)\le  h(x)$ on $D$.
\end{description}
and the $h$-equilibrium measure satisfies
\[ \Vert \gamma \Vert^2=W_h^{-1}(K)=C^{\cal A}_h(K)=\gamma(K). \]
Hence, the $h$-equilibrium potential coincides with the smoothed reduction of $h$ on $K$:
\begin{equation}\label{hpotential=hreduction}
U_{\cal A}^{\gamma}(x)= \hat{R}^{K}_h(x)=\int\limits_{K} \frac{G_{\cal A}(x,y)}{h(y)}\gamma(dy) 
\end{equation}

The following theorem clarifies the variational nature of the $h$-equilibrium measure
and the $h$-equilibrium potential of compact subset $K\subset D$:
\begin{theorem}\label{equilibriummeasure}
Assume $C^{\cal A}_h(K)>0$. Then the $h$-equilibrium measure $\gamma$ may be viewed as the unique solution of the following variational problems:
\begin{description}
\item{(i)\ } In the subset of measures $\mu \in {\cal M}_K$ satisfying the condition
\[ \sup\limits_{S(\mu)}\frac{U_{\cal A}^\mu(x)}{h(x)}=1, \]
find the measure $\gamma$ with $\gamma(K)=\max \mu(K)$.
\item{(ii)\ }In the subset of measures $\mu \in {\cal M}_K$ normalized by the condition
\[ \mu(K)=C_h^{\cal A}(K) \]
find the measure $\gamma$, for which
\[ \sup\limits_{S(\gamma)}\frac{U_{\cal A}^\gamma(x)}{h(x)}=\min \sup\limits_{S(\mu)}\frac{U_{\cal A}^\mu(x)}{h(x)}. \]
\item{(iii)\ } Among all the measures $\mu \in {\cal M}_K$, find the measure $\gamma$, for which
\[ \Vert \gamma \Vert^2-2\gamma(K)=\min \{\Vert \mu \Vert^2-2\mu(K)\} \]
\end{description}
\end{theorem}

The proof of Theorem~\ref{equilibriummeasure} coincides with the proof given in the case
$h=1$ (see Theorem 2.3 in \cite{Landkov}). Note that since the maximum in the variational problem
$(i)$ is equal to $C^{\cal A}_h(K)$, one can give an equivalent defintion of $C^{\cal A}_h(K)$:
\[ C^{\cal A}_h(K)=\max\{\mu(K): \mu \in {\cal M}_K, \ \ U_{\cal A}^\mu(x)\le h(x) \quad\text{for}\  \ x\in S(\mu) \}  \]
Again, this definition coincides with the classical equivalent definition of the Greenian
capacity in the case ${\cal A}=\Delta, h\equiv 1$.

Yet another equivalent definition of the $h$-Capacity is given in Definition 1.4 of Section \ref{mainresult}.
Since $ {\cal G}$ is a closed and convex subset of the Hilber space $H_0^{1,2}(D)$, there exists a unique $u \in {\cal G}$ such that
\begin{equation}\label{capacitarypotential}
 C^{\cal A}_h(K)={\cal L}_{\cal A}(u)
\end{equation}

 {\bf Definition 1.8.}\ The function $u$ which satisfies (\ref{capacitarypotential}) is called the $h$-capacitary potential 
of the compact $K\subset D$ (with respect to $D$ and ${\cal A}$).

The following theorem characterizes $h$-capacitary potential, and demonstrates the equivalence of Definitions 1.4 and 1.6.
\begin{theorem}\label{cappotential}
The $h$-capacitary potential $u$ of any compact $K \subset D$ coincides with the $h$-equilibrium potential 
$U_{\cal A}^{\gamma}(x)$ and with the smoothed $h$-reduction $\hat{R}^{K}_h(x)$ according to (\ref{hpotential=hreduction}).
It is a weak solution of (\ref{Eq:W:1:1}) in $D \setminus K$, and has the boundary value $h(x)$ on $\partial K$, and $0$ on $\partial \Sigma$
in the sense of $H^{1,2}$. Moreover, the $h$-capacitary potential is ${\cal A}$-{\it superharmonic} in $\Sigma$, and it is a weak solution of the problem
(\ref{Eq:W:2:1}) with $\mu$ replaced with the Radon measure 
\begin{equation}\label{gamma_h}
\gamma_h(dy)=\frac{\gamma(dy)}{h(\cdot)}, S(\gamma_h) \subseteq K; \quad\text{i.e.} \ \gamma_h(A)=\int_D \frac{\gamma(dy)}{h(y)}, A\subset \Sigma.
\end{equation}
\end{theorem}
Proof. The proof is similar to proofs given in the case $h \equiv 1$ \cite{LSW}. We have
\[ \min \{u; \hat{R}^{K}_h(x) \} \in H_0^{1,2}(\Sigma). \]
Moreover, since
\[ {\cal L}_{\cal A}( \min \{u; \hat{R}^{K}_h(x) \} ) \leq  {\cal L}_{\cal A}(u) \]
it follows that 
\[ u(x) = h(x), \quad\text{on} \ K  \quad\text{in the sense of} \ H_0^{1,2}(D). \]
With the standard variational technique one can prove that the minimizer $u$ satisfies the inequality (\ref{Eq:W:1:4}),
and hence it is ${\cal A}$-{\it superharmonic} in $\Sigma$. In particular, by choosing the test function $\phi$ with compact support in $\Sigma\setminus K$, it follows that $u$ satisfies (\ref{Eq:W:1:3}). Hence, $u$ is a weak solution of (\ref{Eq:W:1:1}) in $D \setminus K$, and has the boundary value $h(x)$ on $\partial K$, and $0$ on $\partial \Sigma$
in the sense of $H^{1,2}$. The last assertion follows from the Lemma~\ref{2} - (4).

From the Definition 1.6, \eqref{Eq1:6}, \eqref{Eq:W:2:3} it follows that the $h$-Capacities of the compact set $K\subset D$ 
defined for any uniformly elliptic operators ${\cal A}$ and $\ov{{\cal A}}$ satisfy the estimation
\begin{equation}\label{h-capacitiesElliptic}
\lambda^{-2}C^2 C^{\ov{\cal A}}_{\ov{h}}(K) \le C^{\cal A}_h(K) \le \lambda^{2}C^{-2} C^{\ov{\cal A}}_{\ov{h}}(K)
\end{equation}

Since the $h$-capacitary potential is ${\cal A}$-harmonic in $D\setminus K$, and $\zeta$ is an isolated boundary point of $D\setminus K$, from \eqref{hpotential=hreduction} it follows that
\begin{equation}\label{hCapacity=reductionatzeta}
C^{\cal A}_h(K)=\gamma(K)=\lim_{x\to\zeta}\int\limits_{K} \frac{G_{\cal A}(x,y)}{h(y)}\gamma(dy)= \hat{R}^{K}_h(\zeta)
\end{equation}

The following lemma expresses the fact that $h$-capacity is a topological precapacity:
\begin{lemma}\label{precapacity}
As a set function defined in the class of compact subsets of $D$, $C_h^{{\cal A}}(\cdot)$ is a topological precapacity. i.e.
\begin{description}
\item{(1)\ } $C_h^{{\cal A}}(\cdot)$ is strongly subadditive:
\begin{description}
\item{(a)\ } $C_h^{{\cal A}}(K_1) \leq C_h^{{\cal A}}(K_2)$ if $K_1 \subset K_2$;
\item{(b)\ } $C_h^{{\cal A}}(K_1\cup K_2)+C_h^{{\cal A}}(K_1\cap K_2) \leq C_h^{{\cal A}}(K_2)+C_h^{{\cal A}}(K_2)$;
\end{description}
\item{(2)\ } If $K_n$ is a monotone sequence of compact sets with compact limit $K$, then
\[  \lim_{n\to +\infty} C_h^{\cal A}(K_n)=C_h^{\cal A}(K) \]
\end{description}
\end{lemma}
The proof of this lemma is similar to the proofs given in the case $h\equiv 1$ \cite{Doob, HKM}. One can exploit the properties 
of reductions and use the relation \eqref{hCapacity=reductionatzeta} \cite{Doob}, or pursue direct proof based on the definition of $h$-Capacity.

As a topological precapacity, $C_h^{{\cal A}}$ has an extension to countably strongly subadditive Choquet capacity $C_h^{{\cal A}}(\cdot)$ relative to the class of compact subsets of $D$. Recall that the set function $C: 2^D \to [0,+\infty]$ is called a Choquet capacity if it satisfies the conditions $(1a),(2)$ of the Lemma \ref{precapacity}. Extend $C_h^{\cal A}(\cdot)$ to $2^D$ as follows:
\[ C_h^{{\cal A}}(A)=\sup \{ C_h^{{\cal A}}(K): K \quad\text{is compact}, \ K\subset A\}, \ A\subset D \quad\text{is open} \]
\[ C_h^{{\cal A}}(A)=\inf \{ C_h^{{\cal A}}(B): B \quad\text{is open}, \ A\subset B\}, \ A\subset D \quad\text{is arbitrary} \]
The celebrated Choquet capacitability theorem states:
\begin{theorem}\label{Choquet}\cite{Choquet} The set function $C_h^{\cal A}(\cdot): 2^D \to [0,+\infty]$ is a countably strongly subadditive Choquet capacity
defined on all subsets of $D$. All the Borel (and even analytic) subsets $E \subset D$ ara capacitable, i.e.
\begin{equation}\label{capacitable}
C_h^{{\cal A}}(E)=\sup \{ C_h^{{\cal A}}(K): K\quad\text{is compact}, \ K\subset E\}
\end{equation}
\end{theorem}
Consider some examples:

{\it Example 1.} For arbitratry $a>b>0$ consider the following sets: 
\[ h_a=\{x\in D: h(x)=a\}, h_{a,b}=\{x\in D: b \le h(x) \le a\}, h^b=\{x\in D: b \le h(x)\} \]
We have
\begin{equation}\label{hballpotential}
 \hat{R}^{h_a}_h(x)=
\left\{
\begin{array}{l}
a, \quad\text{if} \ \ h(x) \ge a,\\
h(x), \quad\text{if} \ \ 0<h(x)<a.
\end{array}\right.
\end{equation}
By using the formula (\ref{hCapacity=reductionatzeta}) we have
\[ \hat{R}^{h_a}_h(\zeta)=a=\gamma(h_a)=C^{\cal A}_h(h_a) \]
We also have
\[  \hat{R}^{h_{a,b}}_h \equiv \hat{R}^{h_a}_h \]
and hence,
\[ C^{\cal A}_h(h_{a,b})=C^{\cal A}_h(h_a)=a \]
In particular, this demonstrates that $h$-capacity is not an additive set function. Since
$h_{n,b}\uparrow h^b$ as $n\uparrow +\infty$, from (2) of Lemma 1.2 and \eqref{capacitable} it follows that
\[ C^{\cal A}_h(h^b)=+\infty \]

{\it Example 2.} Let $y \in D, y\neq \zeta$ is fixed and real numbers $a,b$ satisfy: $0< b \le a <h(y)$.
Consider a set
\[ J_{a,b}=\{x\in D: b \leq G_{\cal A}(x,y) \leq a\}, \  J_a=\{x \in \Sigma: G_{\cal A}(x,y) \geq a\} \]
We have $\zeta \in int J_a$ and 
\begin{equation}\label{J_a_bpotential}
 \hat{R}^{J_{a,b}}_h(x)=\hat{R}^{\partial J_a}_h(x)=
\left\{
\begin{array}{l}
H_{h}^{J_a}(x), \quad\text{if} \ \ x\in J_a,\\
h(x), \quad\text{if} \ \ G_{\cal A}(x,y) \leq a.
\end{array}\right.
\end{equation}
where $H_{h}^{J_a}$ is a solution of the DP in $J_a$ under the boundary function $h |_{\partial J_a}$.
By using the formula (\ref{hpotential=hreduction}) we also have
\begin{equation}\label{J_a_bpotential-1}
\hat{R}^{J_{a,b}}_h(x)=\hat{R}^{\partial J_a}_h(x)= \int\limits_{\partial J_a} \frac{G_{\cal A}(x,z)}{h(z)}\gamma(dz)
\end{equation}
where $\gamma$ is $h$-capacitary measure for both $\partial J_a$ and $J_{a,b}$. Clearly, we have
\[ S(\gamma) \subseteq \partial J_a = \partial J_a \cap J_{a,b} \]
By choosing $x=\zeta$, from (\ref{J_a_bpotential}), (\ref{J_a_bpotential-1}) it follows
\[ H_{h}^{J_a}(\zeta)=\gamma(\partial J_a)=C^{\cal A}_h(\partial J_a)=C^{\cal A}_h(J_{a,b}) \]
By the maximum principle the following estimation follows:
\begin{equation}\label{h-capacity-est-1}
\min_{\partial J_a} h \leq C^{\cal A}_h(\partial J_a), C^{\cal A}_h(J_{a,b}) \leq \max_{\partial J_a} h
\end{equation}
By choosing $x=y$, from (\ref{J_a_bpotential}), (\ref{J_a_bpotential-1}) it follows
\[ H_{h}^{J_a}(y)=a \int_{\partial J_a} \frac{\gamma(dz)}{h(z)} \]
By using maximum principle again the following alternative estimation follows:
\begin{equation}\label{h-capacity-est-2}
a^{-1}(\min_{\partial J_a} h)^2 \leq C^{\cal A}_h(\partial J_a), C^{\cal A}_h(J_{a,b}) \leq a^{-1}(\max_{\partial J_a} h)^2
\end{equation}

{\it Example 3.} Let $y \in D, y\neq \zeta$ is fixed, and the real number $a$ satisfies: $a > h(y)$.
We have $\zeta \in D \setminus J_a$ and 
\begin{equation}\label{J_a_potential-1}
 \hat{R}^{J_{a}}_h(x)=\hat{R}^{\partial J_a}_h(x)=
\left\{
\begin{array}{l}
H_{h}^{\Sigma \setminus J_a}(x), \quad\text{if} \ \ x\in D \setminus J_a,\\
h(x), \quad\text{if} \ \ J_a.
\end{array}\right.
\end{equation}
where $H_{h}^{\Sigma \setminus J_a}$ is a solution of the DP in $D \setminus J_a$ under the boundary function $h |_{\partial J_a}$.
By using the formula (\ref{hpotential=hreduction}) we also have
\begin{equation}\label{J_a_bpotential-2}
\hat{R}^{J_{a}}_h(x)=\hat{R}^{\partial J_a}_h(x)= \int\limits_{\partial J_a} \frac{G_{\cal A}(x,z)}{h(z)}\gamma(dz)
\end{equation}
where $\gamma$ is the $h$-capacitary measure for both $\partial J_a$ and $J_{a}$, and clearly $S(\gamma) \subseteq \partial J_a$.
By choosing $x=\zeta$, from (\ref{J_a_potential-1}), (\ref{J_a_bpotential-2}) it follows
\[ H_{h}^{\Sigma \setminus J_a}(\zeta)=\gamma(\partial J_a)=C^{\cal A}_h(\partial J_a)=C^{\cal A}_h(J_{a}) \]
By the maximum principle, the following estimation follows:
\begin{equation}\label{h-capacity-est-3}
\min_{\partial J_a} h \leq C^{\cal A}_h(\partial J_a), C^{\cal A}_h(J_{a}) \leq \max_{\partial J_a} h
\end{equation}
By choosing $x=y$, from (\ref{J_a_potential-1}), (\ref{J_a_bpotential-2}) it follows
\[\hat{R}^{J_{a}}_h(y)=\hat{R}^{\partial J_a}_h(y)= h(y)=a \int_{\partial J_a} \frac{\gamma(dz)}{h(z)} \]
By using the maximum principle again, the following alternative estimation follows:
\begin{equation}\label{h-capacity-est-4}
a^{-1}h(y)\min_{\partial J_a} h \leq C^{\cal A}_h(\partial J_a), C^{\cal A}_h(J_{a}) \leq a^{-1}h(y) \max_{\partial J_a} h
\end{equation}

In the next two lemmas we characterize ${\cal A}$-minimal fine topology near the point $\zeta$. Lemma \ref{minimalthinness1} is known for the classical case ${\cal A}=\Delta$ (\cite{Doob}). In particular, Lemma \ref{minimalthinness1} implies the equivalence of Definitions 1.3 and 1.3$^\prime$. 

\begin{lemma}\label{minimalthinness1}
Let $\Omega\subset D$ be an open set and $B$ denotes a Euclidean deleted neigborhood of $\zeta$. Then the following are true:
\begin{description}
\item{(1)\ } Zero-One law
\begin{equation}\label{zero-one-law}
\lim_{B\downarrow \zeta} \hat{R}^{\Omega^c\cap B}_h(\xi)=ch(\xi), \ \xi\in D; \ \quad\text{and either} \ c=0 \ \quad\text{or} \ \ c=1;
\end{equation}
\item{(2)\ } It is $c=0$ (or $c=1$) in \eqref{zero-one-law} if and only if for some (and equivalently for all) $B$
\begin{equation}\label{zero-one-law1}
\hat{R}^{\Omega^c\cap B}_h \not\equiv h \ (\quad\text{or} \  \hat{R}^{\Omega^c\cap B}_h\equiv h), \ \quad\text{on} \ D; 
\end{equation}
\end{description}
\end{lemma}
Proof. $(1)$. Since
\[ 0\leq \hat{R}^{\Omega^c\cap B_1}_h(\xi) \leq \hat{R}^{\Omega^c\cap B_2}_h(\xi), \ \xi \in D; \ \quad\text{if} \ B_1 \subset B_2 \]
the limit
\begin{equation}\label{Aharmoniclimit}
\lim_{B\downarrow \zeta} \hat{R}^{\Omega^c\cap B}_h(\xi)=\tilde{h}(\xi), \ \xi\in D 
\end{equation}
exists and as a finite limit of decreasing ${\cal A}$-harmonic functions $\tilde{h}$ is ${\cal A}$-harmonic in $D$. Moreover, we have
\begin{equation}\label{hisminimalharmonic} 
0 \leq \tilde{h} \leq h, \ \quad\text{on} \ D 
\end{equation}
If ${\cal A}=\Delta$, since $h$ is minimal harmonic on $D$, from \eqref{hisminimalharmonic} it follows that
\begin{equation}\label{hisminimalharmonic1} 
 \tilde{h} =c h, \ \quad\text{on} \ D \ \quad\text{with} \ 0\leq c \leq 1 
\end{equation}
To prove \eqref{hisminimalharmonic1} in general, first note that the limit \eqref{Aharmoniclimit} is true at the point $\zeta$ as well, and according to the Fundamental Convergence Theorem for ${\cal A}$-superharminc functions \cite{Doob, HKM}, $\tilde{h}$ is ${\cal A}$-superharmonic in a unit disk $\Sigma$. Hence it satisfies \eqref{Eq:W:1:4}. According to the Schwartz Theorem \cite{Schwartz}
there exists a nonnegative Radon measure $\mu$ on $\Sigma$ such that
\[ \int_\Sigma a_{ij}\tilde{h}_{x_i}\phi_{x_j}d\xi=\int_\Sigma \phi(\xi) \mu(d\xi) \]
for all $\phi \in C_0^\infty(\Sigma)$. Therefore, $\tilde{h}$ is the unique weak solution of the problem \eqref{Eq:W:2:1} given by \eqref{Eq:W:2:2}, and $\tilde{h}\in H^{1,p}_0(\Sigma)$ for all $p<2$. Since $\tilde{h}$ is ${\cal A}$-harmonic on $D$, the measure $\mu$ can be supported by $\{\zeta\}$ only. Finally, from 
\eqref{hisminimalharmonic}, \eqref{hisminimalharmonic1} follows.

To prove the Zero-One Law, fix some deleted neigborhood $B_1$ of $\zeta$. For all deleted neigborhoods $B\subset B_1$ we have
\begin{equation}\label{idempotent}
\hat{R}^{\Omega^c\cap B}_h=\hat{R}^{\Omega^c\cap B}_{\hat{R}^{\Omega^c\cap B}_h} \leq \hat{R}^{\Omega^c\cap B_1}_{\hat{R}^{\Omega^c\cap B}_h},
\end{equation}
where the equality expresses the fact that taking the smoothed reduction is an idempotent operation, and the inequalty follows from the monotonicity of the smoothed reduction with respect to the set of reduction. Passing to the limit as $B\downarrow \zeta$ in \eqref{idempotent}, due to \eqref{Aharmoniclimit}, \eqref{hisminimalharmonic1} we have
\begin{equation}\label{chleqc2h}
ch \leq \hat{R}^{\Omega^c\cap B_1}_{ch}=c\hat{R}^{\Omega^c\cap B_1}_h
\end{equation}
Now passing to limit as $B_1\downarrow \zeta$, from \eqref{chleqc2h} it follows
\[ c\leq c^2  \Longrightarrow c=0 \quad\text{or} \ c=1. \]\\
$(2)$. If $c=0$ in \eqref{zero-one-law}, it follows that for some deleted neigborhood $B_1$ of $\zeta$
\begin{equation}\label{axirinciashirim}
\hat{R}^{\Omega^c\cap B_1}_h \not\equiv h \ \quad\text{on} \ D
\end{equation}
Due to the monotonicity of the smoothed reduction with respect to the set of reduction, \eqref{axirinciashirim} is true for arbitrary deleted neigborhood $B\subset B_1$. Let $B_2$ be an arbitrary deleted neigborhood of $\zeta$ such that $B_1 \subset B_2$. Prove that \eqref{axirinciashirim} is true with $B_1$ replaced by $B_2$. Assuming on the contrary, and by using the subadditivity of a smoothed reduction we have
\begin{equation}\label{axirinciashirim1}
h \equiv \hat{R}^{\Omega^c\cap B_2}_h \leq \hat{R}^{\Omega^c\cap B_1}_h+\hat{R}^{\Omega^c\cap B_3}_h \ \quad\text{on} \ D
\end{equation}
where $B_3=B_2 \setminus B_1$. Let
\[ \hat{R}^{\Omega^c\cap B_i}_h(\xi) = \int_{\overline{\Omega^c\cap B_i}}G_{{\cal A}}(\xi,\eta)\mu_i(d\eta) \]
where measure $\mu_i$ is supported on a compact subset ${\overline{\Omega^c\cap B_i}}$ of the unit disk $\Sigma$. Clearly, $\mu_2=\delta_\zeta$ is a Dirac point mass of $\{\zeta\}$. Measure $\mu_3$ is compactly supported in $D$, and accordingly $\hat{R}^{\Omega^c\cap B_3}_h$ is an ${\cal A}$-potential in $D$. Assume that
\[ \mu_1(\{\zeta\})=\alpha \]
Due to \eqref{axirinciashirim} we have $0\leq \alpha <1$. From \eqref{axirinciashirim1} it follows that
\begin{equation}\label{axirinciashirim2}
(1-\alpha)h  \leq \int_{\overline{\Omega^c\cap B_i}\setminus \{\zeta\}}G_{{\cal A}}(\xi,\eta)\mu^\prime_1(d\eta)+\hat{R}^{\Omega^c\cap B_3}_h \ \quad\text{on} \ D
\end{equation}
where $\mu^\prime_1$ is the restriction of $\mu_1$ to $\overline{\Omega^c\cap B_i}\setminus \{\zeta\}$. Since the support of $\mu^\prime_1$ is in $D$, the first term on the right hand side of \eqref{axirinciashirim2} is an ${\cal A}$-potential in $D$. Hence, the sum on the right hand side of \eqref{axirinciashirim2} is an ${\cal A}$-potential in $D$. Since the greatest ${\cal A}$-harmonic minorant of the ${\cal A}$-potential is zero, it follows that necessarily we must have $\alpha=1$ in \eqref{axirinciashirim2}, which is a contradiction. Hence, the "only if" claim of $(2)$ is proved. To prove the "if" claim assume that for some $B_1$ and $\xi \in D$ we have 
\begin{equation}\label{axirinciashirim3}
\hat{R}^{\Omega^c\cap B_1}_h(\xi) < h(\xi) 
\end{equation}
For all $B\subset B_1$ we have
\begin{equation}\label{axirinciashirim4}
\hat{R}^{\Omega^c\cap B}_h(\xi) \leq \hat{R}^{\Omega^c\cap B_1}_h(\xi) < h(\xi) 
\end{equation}
Passing to limit along $B \downarrow \zeta$ in \eqref{axirinciashirim4} from the Zero-One Law it follows that $c<1$, and hence $c=0$. Lemma is proved.\\

{\bf Remark:} Note that $\Omega^c$ is ${\cal A}$-minimally thin in ${\cal A}$-minimal fine topology if and only if $c=0$ in Lemma \ref{minimalthinness1}.\\

The following lemma demonstrates that ${\cal A}$-minimal thinness at $\zeta$, and the accordingly minimal-fine topology neigborhood base of 
$\zeta$ is independent of the elliptic operator ${\cal A}$ with a uniform ellipticity constant $\lambda$.
\begin{lemma}\label{minimalthinness}
If $E\subset D$ is ${\cal A}$-minimally thin at $\zeta$, then it is $\overline{{\cal A}}$-minimally thin at $\zeta$ for any elliptic operator $\overline{{\cal A}}$
with uniform ellipticity constant $\lambda$.
\end{lemma}
Proof. Without loss of generality, assume that $E$ is a compact subset of the closed disk $\Sigma_\delta=\{\xi: |\xi|\le 1-\delta\}$ and the constant $C\geq 1$ is fixed to fulfill \eqref{Eq:W:2:3} from Lemma \ref{2} concerning operators ${\cal A}$ and $\overline{{\cal A}}$. Let
\[ E^n=E\cap \{\xi: h(\xi) \geq n\},   \]
and
\[ \hat{R}_h^{E^n}(\xi)=\int_{E^n}\frac{G_{{\cal A}}(\xi,\eta)}{h(\eta)} \gamma_n^h(d\eta), \]
where $\gamma_n^h$ is an $h$-capacitary measure. Since $E$ is ${\cal A}$-minimally thin at $\zeta$, \eqref{Eq:W:1:10a} is fulfilled.
Consider $\overline{{\cal A}}$-potential
\[ \Gamma(\xi)=\int_{E^n}\frac{G_{\overline{{\cal A}}}(\xi,\eta)}{\overline{h}(\eta)} \gamma_n^h(d\eta) \]
From \eqref{Eq:W:2:3} and properties of reductions it follows
\[  \Gamma(\xi) \geq C^{-2}  \hat{R}_h^{E^n}(\xi) = C^{-2} h(\xi) \geq C^{-3} \overline{h}(\xi), \ \xi \in E^n \]
in the sense of $H_0^{1,2}(\Sigma)$. It follows that
\[  \hat{\overline{R}}_{\overline{h}}^{E^n}(\xi) \leq C^3 \Gamma(\xi) \leq C^5  \hat{R}_h^{E^n}(\xi), \ \xi \in D \]
where  $\hat{\overline{R}}$ denotes smoothed reduction with respect to operator $\overline{{\cal A}}$.
Passing to limit as $n\uparrow +\infty$ it follows that \eqref{Eq:W:1:10a} is fulfilled for $\hat{\overline{R}}_{\overline{h}}^{E\cap B}$
for a particular choice of decreasing sequence of neigborhoods $B=\{\xi: h(\xi) \geq n\} \downarrow \zeta$. Lemma is proved.

\subsection{Historical Comments}\label{E:1:2}
It would be convenient to make some remarks concerning the Dirichlet problem for uniformly elliptic equations.  The solvability, in some generalized sense, of the classical DP in a bounded open set $E \subset \rr^{N}$, with prescribed data on $\partial E$, is realized within the class of resolutive boundary functions, identified by Perron's method and its Wiener \cite{Wiener1,Wiener2} and Brelot \cite{Brelot} refinements.
Such a method is referred to as the PWB method, and the corresponding solutions are PWB solutions. 

Wiener, in his pioneering works \cite{Wiener1,Wiener2}, proved a necessary and sufficient condition for the finite boundary point $x_o \in \partial E$ to be regular in terms of the ``thinness" of the complementary set in the neighborhood of $x_o$. Cartan pointed out that the thinness could be characterized by means of
fine topology -- the coarsest topology of $\rr^N$ which makes every superharmonic function continuous.
In fact, a finite boundary point is irregular if and only if $E$ is a deleted neigborhood 
of $x_o$ in fine topology \cite{Doob}. 

De Giorgi \cite{DeGiorgi} and Nash \cite{Nash} almost simultaneously 
proved that any local weak solution of (\ref{Eq:W:1:1}) is locally H\"older 
continuous. Moser \cite{Moser} gave a simpler proof of this fact, as well as the
Harnack inequality. In \cite{LSW} it is proved that the Wiener test for the regularity
of finite boundary points with respect to elliptic operator (\ref{Eq:W:1:1})
coincides with the classical Wiener test for the boundary regularity of harmonic functions.
Hence, the fine-topological neigborhood base of the finite boundary point is 
independent of elliptic operator (\ref{Eq:W:1:1}). The Wiener test for the regularity of finite boundary points for linear degenerate elliptic equations is proved in \cite{FJK}. 
The Wiener test for the regularity of finite boundary points for quasilinear 
elliptic equations was settled due to \cite{Mazya3, gariepy, martio, kilpi}.
Nonlinear potential theory was developed along the same lines as classical potential theory
for the Laplace operator, for which we refer to monographs \cite{HKM, Ziemer}.

In \cite{Abdulla1,Abdulla2,Abdulla3,Abdulla4} a new notion of regularity 
of the boundary point $\infty$ for elliptic and parabolic equations was introduced.
The regularity of $\infty$, for an arbitrary open set $\Omega$, charecterizes the uniqueness of bounded solutions in 
$\Omega$; continuity of solutions at $\infty$; thinness in fine topology at $\infty$; nature of the harmonic measure of $\{\infty\}$; and asymptotic laws for the associated Markov processes at $\infty$.
In \cite{Abdulla4}, the Wiener test for the regularity of $\infty$ was proved for the elliptic PDE (\ref{Eq:W:1:1})
when the space dimension $N\geq 3$, the test being the same as for the Laplace operator \cite{Abdulla1}.

However, the developed approach is not applicable to elliptic PDEs in the 2D case. In particular, the Dirichlet problem in any Greenian open subset of $\rr^2$ always has a unique bounded solution, independetly of the geometry of the boundary at $\infty$. The goal of this paper is to analyze the 2D case. We introduce the notion of $log$-regularity of the boundary point $\zeta$ (finite or $\infty$), for an arbitrary Greenian open set, and prove the Wiener test for the $log$-regularity of $\zeta$. The $log$-regularity of $\zeta$ characterizes the removabilty of the logarithmic singularity at $\zeta$; the uniqueness of solutions in $O(log|\cdot - \zeta|)$ class; the thinness in minimal fine topology near $\zeta$; the nature of $log$-harmonic measure of $\zeta$; and the asymptotic laws for the $log$-Brownian motion.

\section{Proof of Theorem 1.1}\label{proof of the main theorem}
${\bf (1)\Leftrightarrow(2)}$: Assume that $\zeta$ is $log$-{\it regular} and let $u_1$ and $u_2$ be 
two bounded solutions of $h$-DP. Then $v=u_1-u_2$
is a bounded solution of $h$-DP with zero boundary data on $\partial \Omega \setminus \{\zeta\}$. 
Since a generalized solution is order preserving (\cite{Doob, HKM}), we have

\[ |v|\le \ ^hH_{M\cdot1_{\{\zeta\}}}^{\bf \Omega} \equiv M \ ^hH_{1_{\{\zeta\}}}^{\bf \Omega} \equiv 0, \quad\text{with}\> M=\sup|v|. \]
On the other hand, if $\zeta$ is irregular, then for an arbitrary real number $r$, $r \cdot \ ^hH_{1_{\{\zeta\}}}^{\bf \Omega}$ is a generalized solution of the $h$-DP with zero boundary data on $\partial\Omega \setminus \{\zeta\}$.\\

${\bf (2)\Rightarrow(3)}$: Assume that (3) is not satisfied. That is to say, there is a bounded function $f$ such that
for some generalized solution $^hH_f^{\Omega}$, 
one of the inequalities in (\ref{Eq:W:1:10}) is violated. In this case, by choosing a number $\ov{f}$ satisfying 

\[ f_{*}\equiv \liminf_{z\to\infty} f(z) \le \ov{f} \le \limsup_{z\to\infty} f(z)\equiv f^{*}, \]
and by extending $f(\infty)=\ov{f}$, we can always construct a generalized solution $^hH_f^{\Omega}$ which satisfies (\ref{Eq:W:1:10}) . Indeed, since a generalized solution is order preserving, we clearly have
\begin{equation}\label{Eq:W:3:1}
|^hH_f^{\Omega}| \le M \equiv \sup |f|.
\end{equation}
Then for an arbitrary $\epsilon > 0$ we choose $R>0$ such that 
\begin{equation}\label{Eq:W:3:2}
f_* - \epsilon \leq f \leq f^* + \epsilon, \quad\text{ on }\> \partial \Omega \cap \Sigma',
\end{equation}
where $\Sigma' \equiv \{ z \mid h(z) > R \}$ is an open neigborhood of $\zeta$ in a unit disk.
We have
\begin{equation}\label{Eq:W:3:2a}
\hat{R}^{\partial \Sigma'}_h(z)\equiv h(R) = \ ^1H_h^{\Sigma'}(z), \ ^hH_{1_{\partial\Sigma'}}^{\Sigma'}(z)=\frac{^1H_h^{\Sigma'}(z)}{h(z)}=\frac{h(R)}{h(z)},  \quad\text{for}\ z\in \Sigma'
\end{equation} 
Since the generalized solution is order preserving, we also have
\begin{equation}\label{Eq:W:3:3}
f_* - \epsilon - 2M \ ^hH_{1_{\partial\Sigma'}}^{\Sigma'} \leq \ ^hH_f^{\Omega} \leq  f^* + \epsilon + 2M \ ^hH_{1_{\partial\Sigma'}}^{\Sigma'} \quad\text{ on }\> \Omega \cap \Sigma'.
\end{equation}
This follows from the fact that (\ref{Eq:W:3:3}) is satisfied on $\partial (\Omega\cap \Sigma')$. 
Passing to limit, first as $z\rightarrow \zeta$, and then as $\epsilon \downarrow 0$, from (\ref{Eq:W:3:3}) and (\ref{Eq:W:3:2a}) it follows that the constructed generalized solution satisfies  (\ref{Eq:W:1:10}). Contradiction with uniqueness.\\  

${\bf (3)\Rightarrow(2)}$:  Let $u_1$ and $u_2$ be two bounded solutions of the $h$-DP. Their
difference is a generalized solution with a zero boundary function, and accordingly,
it vanishes identically in view of (\ref{Eq:W:1:10}). On the other hand, if at least for one $f$, (\ref{Eq:W:1:10})
is violated, then from the relation $(2)\Rightarrow(3)$ for direct assertions, it follows that there are
at least two solutions of the $h$-DP. This implies that the $h$-DP with zero boundary data on $\pl \Omega$ has a
non-trivial solution. That, in turn, implies that it must have infinitely many solutions. \\ 

${\bf (1)\Leftrightarrow(4)}$: Assume that $\zeta$ is $log$-regular for $\Omega$. From (\ref{Eq:W:1:8e}) it follows that
\begin{equation}\label{Eq:W:3:2b}
1\equiv \ ^hH_{1_{D\cap\partial\Omega}}^\Omega(\xi)=\mu_\Omega^h(\xi,1_D)=\frac{\mu_\Omega^1(\xi,h1_D)}{h(\xi)}=\frac{^1H_{h1_D}^{\Omega}(\xi)}{h(\xi)}, \ \ \xi\in\Omega
\end{equation}
Moreover, it is easy to prove that
\begin{equation}\label{Eq:W:3:2c}
^hH_{1_{D\cap\partial\Omega}}^\Omega(\xi)=\frac{^1H_{h1_D}^{\Omega}(\xi)}{h(\xi)}=\frac{\hat{R}^{D-\Omega}_h(\xi)}{h(\xi)}= \ ^h\hat{R}^{D-\Omega}_1(\xi), \ \ \xi\in\Omega
\end{equation}
The relation between Perron solution and reduction expressed in (\ref{Eq:W:3:2c}) is proved in \cite{Doob} (page 122) for the Laplace equation.  The proof applies to general elliptic PDEs, including (\ref{Eq:W:1:1}) without any change (see \cite{HKM}, page 187). From (\ref{Eq:W:3:2b}),(\ref{Eq:W:3:2c}) it follows that
\[ \hat{R}^{D-\Omega}_h \equiv h \ \ \xLongrightarrow{Lemma \ \ref{minimalthinness1}}\quad \  \zeta \in (\Omega^c\cap D)^{mf}  \]
On the contrary, let $\zeta$ be $log$-irregular for $\Omega$. From (\ref{Eq:W:1:8e}) it follows that
\begin{equation}\label{Eq:W:3:2d}
1\not\equiv \ ^hH_{1_{D\cap\partial\Omega}}^\Omega(\xi), \ \ \xi\in\Omega
\end{equation}
From (\ref{Eq:W:3:2d}),(\ref{Eq:W:3:2c}) it follows that
\[ \hat{R}^{D-\Omega}_h(\xi)\not\equiv h(\xi), \ \xi\in\Omega \ \ \xLongrightarrow{Lemma \ \ref{minimalthinness1}}\quad \ \zeta \not\in (\Omega^c\cap D)^{mf}  \]

${\bf (4)\Rightarrow(5)}$: Note that the equivalence of the convergence of two series in \eqref{Eq:W:1:12} follows from \eqref{hCapacity=reductionatzeta}. Assume that the series (\ref{Eq:W:1:12}) converges. Prove that $\Omega^c$ is ${\cal A}$-minimally thin at $\zeta$. 
Let us fix  $\xi^*\in D$ and choose $N_0$ so large that 
\[ \xi^* \in D \cap \{x:  h(x)\leq a^{N_0-1}  \}  \]
Let $E^N=\sum_{n=N}^{+\infty}E_n$. Since $\hat{R}^{\cdot}_h(\xi^*)$ is a countably subadditive set function, for all $N\ge N_0$ we have
\begin{equation}\label{Eq:W:3:3c}
\hat{R}^{E^N}_h(\xi^*)=\hat{R}^{\bigcup_{n=N}^{+\infty} E_n}_h(\xi^*) \le \sum_{n=N}^{+\infty}  \hat{R}^{E_n}_h(\xi^*) 
=\sum_{n=N}^{+\infty}\int\limits_{E_n} \frac{G_{{\cal A}}(\xi^*,y)}{h(y)}\gamma_{E_n}(dy),
\end{equation}
where $\gamma_{E_n}$ is the $h$-equilibrium measure of $E_n$. By using (\ref{Eq:W:2:3c}), from (\ref{Eq:W:3:3c}) it follows that
\begin{equation}\label{Eq:W:3:3a}
\hat{R}^{E^N}_h(\xi^*)\le \sup_{y\in E^N} G_{\cal A}(\xi^*,y) \sum_{n=N}^{+\infty} a^{-n}C_h^{\cal A}(E_n) \le C log\frac{2}{d(N_0)} \sum_{n=N}^{+\infty} a^{-n}C_h^{\cal A}(E_n)
\end{equation}
where $d(N_0)=dist(\Gamma(N_0), \Gamma(N_0-1), \Gamma(N_0)= \{\xi: h(\xi)=N_0 \}$. Since the series (\ref{Eq:W:1:12}) converges by choosing $N$ sufficiently large from (\ref{Eq:W:3:3a}) it follows that
\begin{equation}\label{Eq:W:3:3b}
\lim_{N\uparrow +\infty}\hat{R}^{E^N}_h(\xi^*)=0
\end{equation}
Hence, from the Lemma \ref{minimalthinness1} it follows that $\Omega^c$  is ${\cal A}$-minimally thin at $\zeta$.

${\bf (5)\Rightarrow(4)}$ Assume that $\Omega^c$ is ${\cal A}$-minimally thin at $\zeta$. Prove that both series in \eqref{Eq:W:1:12} converge. Let the constant $C\geq 1$ is fixed to satisfy \eqref{Eq:W:2:3a} in a compact subset $\overline{E}$ of $\Sigma$, where
\[ E:=\Omega^c \cap \{\xi: h(\xi) \geq a \} = \bigcup_{n=1}^{+\infty}E_n  \]
Next, we choose smallest integer $p\geq 2$ which satisfy
\begin{equation}\label{subseries}
C^2 < a^{p-1}
\end{equation}
It is sufficient to prove that any of the following subseries converge:
\begin{equation}\label{subseries1}
\sum_n a^{-np+i} \hat{R}^{E_{np+i}}_h(\zeta), \ i=0,1,...,p-1;
\end{equation}
We are going to prove the convergence of the subseries \eqref{subseries1} with $i=0$. The proof of the convergence of the subseries \eqref{subseries1} with $i=1,...,p-1$ is the same. 
Let 
\[ E^0= \bigcup_{n=1}^{+\infty}E_{np}  \]
Since $E^0\subset \Omega^c$, it is also ${\cal A}$-minimally thin at $\zeta$. By Lemma \ref{minimalthinness} $E^0$ is minimally thin at $\zeta$. From \cite{Doob}, Chapter XII, Section 13 it follows that there is a positive superharmonic function $v$ such that
\begin{equation}\label{mf-limit}
\lim_{E^0 \ni \eta \to \zeta} \frac{v(\eta)}{\log\frac{1}{|\eta|}}=+\infty,  mf-\lim_{\eta \to \zeta} \frac{v(\eta)}{\log\frac{1}{|\eta|}} = inf_{D} \frac{v(\eta)}{\log\frac{1}{|\eta|}} = M_v(\{\zeta\}) 
\end{equation}
where $M_v(\cdot)$ is the measure supported on $\partial D$, associated with the harmonic component of the Riesz decomposition of $v$. In fact, we can replace $v$ with $v(\eta)+M_v(\{\zeta\}\log|\eta|$ to have \eqref{mf-limit} with $M_v(\{\zeta\})=0$ on the right-hand side. Note that if the intersection of $E^0$ with any neigborhood of $\zeta$ is a polar set, then clearly the subseries \eqref{subseries1} with $i=0$ is convergent. Hence without loss of generality we assume that the intersection of $E^0$ with any neigborhood of $\zeta$ is not a polar set. Let us further replace $v$ with its smoothed reduction $\hat{R}_v^{E^0}$, which is also a positive superharmonic function on $\Sigma$, and satisfies the modified \eqref{mf-limit}, where $M_v(\{\zeta\})=0$, and the limit $+\infty$ holds when $\zeta$ is approached through an $E^0$-less polar set. Moreover, since $E^0$ is a relatively compact subset of $\Sigma$, $\hat{R}_v^{E^0}$
is a potential with associated Riesz measure $\mu$: 
\begin{equation}\label{Rieszpotential}
\hat{R}_v^{E^0}(\xi)= \int_{\overline{E^0}}G(\xi,\eta)\mu(d\eta)
\end{equation}
where $G$ is a Green function \eqref{Greenfunction}.  Let $u$ be a solution of the problem \eqref{Eq:W:2:1} with measure $\mu$. 
According to \cite{LSW}, there exists a unique solution such that $u\in H_0^{1,p}(\Sigma), \forall p<2$ and (see Lemma \ref{2}):
\begin{equation}\label{Representationformula_1}
u(x)=\int_{\overline{E^0}} G_{{\cal A}}(x,y)\mu(dy).
\end{equation}
Due to estimation \eqref{Eq:W:2:3}, $u$ satisfies the modified \eqref{mf-limit}, where $M_v(\{\zeta\})=0$, $-\log|\eta|$ is replaced with $h(\eta)$ and the limit $+\infty$ holds when $\zeta$ is approached through an $E^0$-less polar set. Let $\mu_{pn}$ and $\mu_{pn}'$ be projections of $\mu$ to $E_{pn}$ and $\overline{E^0}-E_{pn}$ respectively. 

Let us introduce the sets $\tilde{E}_{pn}$ which are expansions of the sets $E_{pn}$:
\[ \tilde{E}_{pn}=\Omega^c \cap F_{pn}, \ F_{pn}= \{\xi: C^{-1}a^{pn} \leq \log\frac{1}{|\xi|} \leq Ca^{pn+1} \} \]
Due to \eqref{Eq:W:2:3}, we have $E_{np} \subset  \tilde{E}_{pn}$. It can be verified that due to \eqref{subseries}, the sets  $\tilde{E}_{pn}, n=1,2,...$ are separated. Elementary calculation shows that there exists an integer $n_p$ depending only on $p$ such that for all $n \geq n_p$
\begin{equation}\label{log_estimate}
\log\frac{2}{|\xi - \eta|} \leq C^{\prime} \log\frac{1}{|\xi|}, \ \xi \in  \tilde{E}_{pn}, \ \eta \in \  \bigcup_{k \neq n}\tilde{E}_{pk}
\end{equation}
with $C^{\prime}=C^2a+1$. 
From \eqref{Eq:W:2:3}-\eqref{Eq:W:2:3c}, \eqref{log_estimate} it follows that there exists a constant $C^{\prime \prime}>0$ such that for all $n\geq n_p$
\begin{equation}\label{Greenestimation}
G_{{\cal A}}(\xi,\eta)\leq C^{\prime \prime}h(\xi), \quad\text{for} \ \xi\in E_{pn}, \eta\in \overline{E^0}-E_{pn}
\end{equation}
Hence, for $\xi \in E_{pn}$, we have 
\[ \int_{\overline{E^0}-E_{pn}} \frac{G_{{\cal A}}(\xi, \eta)}{h(\xi)} \mu_{pn}'(d\eta)\leq C^{\prime \prime}\mu_{pn}'(\overline{E^0}-E_{pn})\leq \mu(\overline{E^0})<+\infty. \]
Therefore, for sufficiently large $n$
\[  \int_{E_{pn}}G_{{\cal A}}(\xi,\eta)\mu_{pn}(d\eta) \geq h(\xi)\quad\text{approximately everywhere (a.e.) in} \  E_{pn}  \]
Hence,
\[ \hat{R}_h^{E_{pn}}(\xi) \leq \int_{E_{pn}}G_{{\cal A}}(\xi,\eta)\mu_{pn}(d\eta)  \quad\text{a.e. in} \  E_{pn}  \]
and by the domination principle in all $D$. Passing to limit as $\xi \to \zeta$, we have
\[  \hat{R}_h^{E_{pn}}(\zeta) \leq  \int_{E_{pn}}h(\eta)\mu_{pn}(d\eta) \leq a^{pn+1} \mu_{pn}(E_{pn})  \]
and accordingly,
\[
\sum_n a^{-pn}  \hat{R}_h^{E_{pn}}(\zeta) \leq
\sum_n a^{-pn} a^{pn+1} \mu_{pn}(E_{pn}) 
\leq a \mu(E^0) < +\infty \]

Proof of the equivalence ${\bf (5)\Leftrightarrow(6)}$ is standard.

Equivalence ${\bf (5)\Leftrightarrow(6)\Leftrightarrow(7)}$ follows from the equivalence of ${\bf (5)\Leftrightarrow(4)}$ and Lemma \ref{minimalthinness}.

Finally,  the equivalence of ${\bf (5)\Leftrightarrow(8)}$ is a consequence of the estimation
\[ a^n \hat{R}_1^{E_n} \leq  \hat{R}_h^{E_n} \leq a^{n+1} \hat{R}_1^{E_n} \]
which is true approximately everywhere on $E_n$, and accordingly by the domination principle everywhere on $D$.




\end{document}